\documentclass[12pt]{article}
\usepackage{amsmath}
\usepackage{amsfonts}
\usepackage{amssymb}
\usepackage{theorem}

\title{Actions of metric groups and continuous logic}
\author{A.Ivanov 
\thanks{The research is supported by Polish National Science Centre grant DEC2011/01/B/ST1/01406} 
}

\newtheorem{theorem}{Theorem}[section]
\newtheorem{proposition}[theorem]{Proposition}
\newtheorem{corollary}[theorem]{Corollary}
\newtheorem{lemma}[theorem]{Lemma}
\newtheorem{definition}[theorem]{Definition}
\newtheorem{example}[theorem]{Example}
\newtheorem{remark}[theorem]{Remark}

\begin{document} 
\topmargin = 12pt
\textheight = 630pt 
\footskip = 39pt 

\maketitle

\begin{quote}
{\bf Abstract} 
We study expressive power of continuous logic in classes of 
metric groups defined by properties of their actions. 
For example we consider properties non-{\bf OB}, non-{\bf FH} 
and non-{\bf F}$\mathbb{R}$. 
\end{quote}

\section{Introduction} 
Hereditary properties of basic classes of topological groups 
studied in measurable and geometric group theory 
have deserved more attention of researches in recent 
investigations, see \cite{GriHarp}, \cite{HHM}, \cite{pestov16}. 
This is mainly connected with the tendency of study 
of such notions as amenability or property ${\bf (T)}$ 
of Kazhdan outside the class of locally compact groups, 
see \cite{EvTs}, \cite{hartnickkohl}, \cite{ivarssonkutzschebauch} and \cite{tsankov}.   

From this point of view it is natural to verify 
the behaviour of these classes under logical constructions. 
Moreover this task looks quite attractive because 
some logical constructions, for example ultraproducts, 
have become common in group theory.  

On the other hand typical classes of groups studied in geometric 
group theory are non-axiomatizable. 
For example let us consider the following well-known classes 
of topological groups. 
\begin{definition} \label{def_prop} 
\begin{itemize} 
\item $G \in {\bf FH}$ if any continuous  
affine  isometric  action of $G$ on 
a real Hilbert space has a fixed point;    
\item $G\in {\bf F}\mathbb{R}$ if any continuous 
isometric action of $G$ on a real tree has a fixed point;  
\item $G \in {\bf OB}$ if for any continuous isometric action 
of $G$ on a metric space all orbits are bounded. 
\end{itemize} 
\end{definition} 

We remind the reader that an action of $G$ on a metric space $X$ 
is called continuous if it is continuous as a 2-argument function 
$G\times X \rightarrow X$. 
When the action is isometric this is equivalent to the condition that 
for any $x\in X$ the map $g \rightarrow gx$ is continuous. 
 
Discrete groups of the class {\bf OB} are called {\bf strongly bounded}.  
Y. de Cornulier have proved in \cite{dC} that 
they are contained in {\bf FH} and ${\bf F}\mathbb{R}$. 
Moreover it is also shown in \cite{dC} that for any finite perfect 
group $F$ and an infinite $I$ the power $F^I$ is strongly bounded.  

Since $F^I$ is locally finite, any its countable subgroup 
has cofinality $\omega$, i.e. is the union of a strictly 
increasing $\omega$-chain of proper subgroups.  
Since such groups are outside of 
${\bf FH} \cup {\bf F}\mathbb{R}$, 
any countable elementary subgroup $H$ of $F^I$ 
witnesses the non-axiomatizability of 
{\bf OB}, {\bf FH}, ${\bf F}\mathbb{R}$ and 
non-{\bf OB}, non-{\bf FH}, non-${\bf F}\mathbb{R}$.  

Although this argument is carried out in the discrete 
case, it can be applied in many other situations,   
for example in continuous logic. 
Thus we see that the basic logic constructions 
involving elementary equivalence look foreign  
to the properties defined above. 

On the other hand note that the properties we look at 
are formulated in the language of $G$-actions. 
Thus in order to adapt the situation to the logic 
approach let us consider the following definition 
and the corresponding question after it. 

\begin{definition} \label{def_anal} 
Let ${\cal K}$ be a class of (first-order/continuous) structures. 
We say that ${\cal K}$ is {\bf logically analyzable} 
if there is a family of (first-order/continuously) 
axiomatizable classes 
$\mathcal{K}_{\alpha}$, $\alpha \in I$, 
in expanded languages with possibly new sorts, so that 
\begin{itemize} 
\item for each $\alpha$ reducts of structures of $\mathcal{K}_{\alpha}$ 
to the language of $\mathcal{K}$ belong to  $\mathcal{K}$,  
\item every $G\in \mathcal{K}$ has an expansion $\hat{G}$ 
belonging to one of these classes $\mathcal{K}_{\alpha}$.  
\end{itemize} 
\end{definition} 

We now formulate the main question of the paper. 

\bigskip 

{\bf Analyzability Question.} 
Let ${\cal K}$ be one of the classes 
{\bf OB}, {\bf FH}, {\bf F}$\mathbb{R}$, non-{\bf OB}, 
non-{\bf FH}, non-{\bf F}$\mathbb{R}$ or any other class 
of groups. 
{\em Is $\mathcal{K}$ logically analyzable?}
\bigskip 

Note that we may consider some other kinds of axiomatizability in this question,  
for example the $L_{\omega_1 ,\omega}$-version of continuous logic.

This idea is connected with papers \cite{KM}, \cite{phillips}, 
\cite{sabbagh} and \cite{thomas} where the following 
Ph. Hall's notion is investigated. 
A class of discrete groups $\mathcal{K}$ is called {\bf bountiful} if 
for any pair of infinite groups $G\le H$ with $H\in \mathcal{K}$ 
there is $K\in \mathcal{K}$ such that $G\le K\le H$ and $|G|=|K|$. 
Generalizing some logical observations from 
\cite{KM} it is easy to see that if 
$\mathcal{K}$ is logically analyzable then $\mathcal{K}$ is bountiful 
(see Section \ref{S_1_1} for a precise argument).

When one considers topological groups, the definition of 
bountiful classes should be modified by replacing cardinality 
of groups  by density character, 
i.e. the smallest cardinality of a dense subset of $K$. 
Moreover it is natural to replace the subgroup 
$G$ in the definition by a set. 
As a result we formulate the definition as follows. 
\begin{definition} 
A class of topological groups $\mathcal{K}$ is called  {\bf bountiful} 
if for any infinite group $H\in \mathcal{K}$ and any $C\subseteq H$ 
there is $K\in \mathcal{K}$ such that $C\subseteq  K\le H$ and 
the density character of $C$ 
coincides with the density character of $K$. 
\end{definition} 

Proposition \ref{Bounti}(a) in the final part of this section 
shows that properties {\bf OB}, {\bf FH} and ${\bf F}\mathbb{R}$ 
are not bountiful. 
So, as we have already mentioned, it is easy to see that 
these classes are not logically analyzable. 
On the other hand part \ref{Bounti}(b) of this proposition states that 
properties non-{\bf OB}, non-{\bf FH} and non-${\bf F}\mathbb{R}$
are bountiful, i.e. we may conjecture that 
the classes non-{\bf OB}, non-{\bf FH} and non-{\bf F}$\mathbb{R}$ 
are logically analyzable.  

The main results of our paper confirm this conjecture  
under some uniform continuity assumptions. 
They are presented in Sections 3 and 4. 
As a consequence we obtain bountifulness 
of some uniform versions (i.e. subclasses) of non-{\bf FH} 
and non-{\bf F}$\mathbb{R}$
in the form which is more precise
than the statements of Proposition \ref{Bounti}(b) 
below.  

\paragraph{Novelty of the approach.} 
When we apply logical methods to properties 
involving group actions the basic problem 
which we face is axiomatization of the action. 
Typically unbounded metric spaces are considered in 
continuous logic as many-sorted structures of 
$n$-balls of a fixed point of the space ($n\in \omega$). 
Section 15 of \cite{BYBHU} contains nice examples of 
such structures.

If the action of a bounded metric group $G$ is isometric 
and preserves these balls we may consider 
the action as a sequence of binary operations where  
the first argument corresponds to $G$. 
In such a situation one just fixes 
a sequence of continuity moduli for $G$ 
(for each $n$-ball). 
We will see in Section 2 that this approach  
works well for the negation of property ${\bf (T)}$ 
(non-${\bf (T)}$)  
in the class of locally compact groups. 

The situation dramatically changes 
when the action does not preserve 
$n$-balls. 
For example this happens when we study properties 
{\bf FH}/non-{\bf FH} (or {\bf F}$\mathbb{R}$/non-{\bf F}$\mathbb{R}$), 
where affine actions on Hilbert spaces appear 
(or non-elliptic actions on unbounded trees). 
In Section 3 we present a new approach to such situations. 
Using geometric properties of Hilbert spaces and 
real trees we introduce sequences of ternary predicates 
and show that under some natural assumptions 
on the action, continuity moduli for these predicates can be defined. 
This is the crucial element of the paper. 
It allows us to axiomatize classes of actions which we consider, 
see Theorems \ref{isoR} and \ref{isoH}.  

In Section 4 we slightly simplify the circumstances. 
Replacing non-{\bf OB} by some uniform versions 
of it we arrive at a situation where instead of 
adding new sorts one just adds two continuous predicates 
to the signature. 
This trick can be also applied to non-Roelcke bounded groups, 
non-Roelcke precompact groups and non-${\bf (OB)_k}$-groups 
(see  \cite{rosendalN}).

\paragraph{Uniform continuity.} 
Actions of metric groups which can be analyzed by 
tools of continuous logic must be uniformly continuous 
for each sort appearing in the presentation of the space by metric balls. 
This slightly restricts the field of applications 
of our results. 
On the other hand note that in the case 
of discrete groups we do not lose generality and 
moreover our methods become more powerful. 
In Section 4.B we analyze some other properties  
of discrete groups, for example {\bf FA}.

\bigskip

\begin{remark} \label{Hofmann}
{\em 
We mention paper \cite{HHM} where related questions were studied 
in the case of locally compact groups. 
It is proved in \cite{HHM} that for any locally compact group $G$, 
the entire interval of cardinalities between $\aleph_0$ and $w(G)$, 
the weight of the group, is occupied by the weights of closed 
subgroups of $G$.   
We remind the reader that the weight of a topological space $(X,\tau )$ 
is the smallest cardinality which can be realized as 
the cardinality of a basis of $(X,\tau )$.  
If the group $G$ is metric, the weight of $G$ coincides 
with the density character of $G$. 
This yields the following version of the L\"{o}wenheim-Skolem Theorem 
(see Section \ref{S_1_1}). }

 Let $G$ be a locally compact group which is a continuous structure. 
Then for any cardinality $\kappa < \mathsf{density}(G)$ there is 
a closed subgroup $H<G$ such that $\mathsf{density}(H)=\kappa$ 
and $H$ is an elementary substructure of $G$. 
\end{remark}

In the rest of this introduction we briefly remind 
the reader some preliminaries of continuous logic.  
Since we want to make the paper available for non-logicians 
these preliminaries can look too tedious for specialists. 
On the other hand we inform the reader that all necessary algebraic 
terms will be defined in the introductionary parts of 
corresponding sections.

\subsection{Continuous structures} 
\label{S_1_1}

We fix a countable continuous signature 
$$
L=\{ d,R_1 ,...,R_k ,..., F_1 ,..., F_l ,...\}. 
$$  
Let us recall that a {\em metric $L$-structure} 
is a complete metric space $(M,d)$ with $d$ bounded by 1, 
along with a family of uniformly continuous operations on $M$ 
and a family of predicates $R_i$, i.e. uniformly continuous maps 
from appropriate $M^{k_i}$ to $[0,1]$.   
It is usually assumed that to a predicate symbol $R_i$ 
a continuity modulus $\gamma_i$ is assigned so that when 
$d(x_j ,x'_j ) <\gamma_i (\varepsilon )$ with $1\le j\le k_i$ 
the corresponding predicate of $M$ satisfies 
$$ 
|R_i (x_1 ,...,x_j ,...,x_{k_i}) - R_i (x_1 ,...,x'_j ,...,x_{k_i})| < \varepsilon . 
$$ 
It happens very often that $\gamma_i$ coincides with $id$. 
In this case we do not mention the appropriate modulus. 
We also fix continuity moduli for functional symbols. 
Each classical first-order structure can be considered 
as a complete metric structure with the discrete $\{ 0,1\}$-metric. 

By completeness continuous substructures of a continuous structure are always closed subsets. 

Atomic formulas are the expressions of the form $R_i (t_1 ,...,t_r )$, 
$d(t_1 ,t_2 )$, where $t_i$ are terms (built from functional $L$-symbols). 
In metric structures they can take any value from $[0,1]$.   
{\em Statements} concerning metric structures are usually 
formulated in the form 
$$
\phi = 0 
$$ 
(called an $L$-{\em condition}), where $\phi$ is a {\em formula}, 
i.e. an expression built from 
0,1 and atomic formulas by applications of the following functions: 
$$ 
x/2  \mbox{ , } x\dot- y= \mathsf{max} (x-y,0) \mbox{ , } \mathsf{min}(x ,y )  \mbox{ , } \mathsf{max}(x ,y )
\mbox{ , } |x-y| \mbox{ , } 
$$ 
$$ 
\neg (x) =1-x \mbox{ , } x\dot+ y= \mathsf{min}(x+y, 1) \mbox{ , } \mathsf{sup}_x \mbox{ and } \mathsf{inf}_x . 
$$ 
A {\em theory} is a set of $L$-conditions without free variables 
(here $\mathsf{sup}_x$ and $\mathsf{inf}_x$ play the role of quantifiers). 
   
It is worth noting that any formula is a $\gamma$-uniformly continuous 
function from the appropriate power of $M$ to $[0,1]$, 
where $\gamma$ is the minimum of continuity moduli of $L$-symbols 
appearing in the formula. 

The condition that the metric is bounded by $1$ is not necessary. 
It is often assumed that $d$ is bounded by some rational number $d_0$. 
In this case the (truncated) functions above are appropriately modified.  

We sometimes replace conditions of the form $\phi \dot{-} \varepsilon =0$ 
where $\varepsilon \in [0,d_0]$ by more convenient expressions $\phi \le \varepsilon$. 

A tuple $\bar{a}$ from $M^n$ is {\em algebraic} in $M$ over 
$A$ if there is a compact subset $C\subseteq M^n$ such that 
$\bar{a}\in C$ and the distance predicate $\mathsf{dist}(\bar{x},C)$ 
is definable (in the sense of continuous logic, \cite{BYBHU}) 
in $M$ over $A$. 
Let $\mathsf{acl}(A)$ be the set of all elements algebraic over $A$. 
In continuous logic the concept of algebraicity is 
parallel to that in traditional model theory 
(see Section 10 of \cite{BYBHU}).  

\paragraph{Axiomatizability in continuous logic. } 

When one considers classes axiomatizable in continuous logic 
it is usually assumed that all operations and predicates 
are uniformly continuous with respect to some fixed continuity moduli. 
Suppose that $\mathcal{C}$ is a class of metric $L$-structures. 
Let $Th^c (\mathcal{C})$ be the set of all $L$-conditions 
without free variables which hold in all structures of $\mathcal{C}$. 
It is proved in \cite{BYBHU} (Proposition 5.14 and Remark 5.15) 
that every model of $Th^c (\mathcal{C})$ is elementary equivalent 
to some ultraproduct of structures from $\mathcal{C}$. 

\paragraph{Metric groups.} 
Below we always assume that our metric groups are 
continuous structures with respect 
to bi-invariant metrics (see \cite{BYBHU}). 
This exactly means that $(G,d)$ is 
a complete metric space and $d$ is bi-invariant. 
Note that the continuous logic approach 
takes weaker assumptions on $d$. 
Along with completeness it is only assumed 
that the operations  of a structure 
are uniformly continuous with respect to $d$. 
Thus it is worth noting here that 
any group which is a continuous structure has 
an equivalent bi-invariant metric. 
See \cite{sasza2} for a discussion 
concerning  this observation.

\paragraph{Hilbert spaces in continuous logic.} 
We treat a Hilbert space over $\mathbb{R}$ 
exactly as in Section 15 of \cite{BYBHU}. 
We identify it with a many-sorted metric structure 
$$
(\{ B_n\}_{n\in \omega} ,0,\{ I_{mn} \}_{m<n} ,
\{ \lambda_r \}_{r\in\mathbb{R}}, +,-,\langle \rangle ),
$$
where $B_n$ is the ball of elements of norm $\le n$, 
$I_{mn}: B_m\rightarrow B_n$ is the inclusion map, 
$\lambda_{r}: B_m\rightarrow B_{km}$ is scalar 
multiplication by $r$, with $k$ the unique integer 
satisfying $k\ge 1$ and $k-1 \le |r|<k$; 
futhermore, $+,- : B_n \times B_n \rightarrow B_{2n}$ 
are vector addition and subtraction and 
$\langle \rangle : B_n \rightarrow [-n^2 ,n^2 ]$ 
is the predicate of the inner product. 
The metric on each sort is given by 
$d(x,y) =\sqrt{ \langle x-y, x-y \rangle }$.   
For every operation the continuity modulus is standard.  
For example in the case of $\lambda_r$ this is $\frac{z}{|r|}$. 
Note that in this version of continuous logic 
we do not assume that the diameter of a sort is bounded by 1. 
It can become any natural number.  

Stating existence of infinite approximations of orthonormal bases 
(by a countable family of axioms, see Section 15 of \cite{BYBHU}) 
we assume that our Hilbert spaces are infinite dimensional. 
By \cite{BYBHU} they form the class of models of a complete 
theory which is $\kappa$-categorical for all infinite $\kappa$, 
and admits elimination of quantifiers. 

This approach can be naturally extended to complex Hilbert spaces, 
$$
(\{ B_n\}_{n\in \omega} ,0,\{ I_{mn} \}_{m<n} ,
\{ \lambda_c \}_{c\in\mathbb{C}}, +,-,\langle \rangle_{Re} , \langle \rangle_{Im} ). 
$$
We only extend the family 
$\lambda_{r}: B_m\rightarrow B_{km}$, $r\in \mathbb{R}$, 
to a family $\lambda_{c}: B_m\rightarrow B_{km}$, $c\in \mathbb{C}$, 
of scalar products by $c\in\mathbb{C}$, with $k$ 
the unique integer satisfying $k\ge 1$ and $k-1 \le |c|<k$. 

We also introduce $Re$- and $Im$-parts of the inner product. 

If we remove from the signature of complex Hilbert spaces 
all scalar products by $c\in \mathbb{C}\setminus \mathbb{Q}[i]$, 
we obtain a countable subsignature 
$$
(\{ B_n\}_{n\in \omega} ,0,\{ I_{mn} \}_{m<n} ,
\{ \lambda_c \}_{c\in\mathbb{Q}[i]}, +,-,\langle \rangle_{Re} , \langle \rangle_{Im} ),
$$
which is {\em dense} in the original one: \\ 
- if we present $c\in \mathbb{C}$ by a sequence $\{ q_i \}$ 
from $\mathbb{Q}[i]$ converging to $c$, 
then the choice of the continuity moduli of 
the restricted signature still guarantees that 
in any sort $B_n$ the functions $\lambda_{q_i}$ 
form a sequence which converges to  $\lambda_c$ 
with respect to the metric 
$$ 
\mathsf{sup}_{x\in B_n} \{ |f^M (x) - g^M (x)| : M \mbox{ is a model of the theory of Hilbert spaces } \}.  
$$  

\paragraph{L\"{o}wenheim-Skolem theorem and bountifulness.} 
The following theorem is one of the main tools of this paper. 
\bigskip 
 
{\bf L\"{o}wenheim-Skolem Theorem.} (\cite{BYBHU}, Proposition 7.3) 
{\em Let $\kappa$ be an infinite cardinal number and assume 
$|L|\le \kappa$. 
Let $M$ be an $L$-structure and suppose $A\subset M$ has 
density $\le\kappa$. 
Then there exists a substructure $N\subseteq M$ containing $A$ such that 
$\mathsf{density}(N) \le\kappa$ and $N$ is an elementary substructure of $M$, i.e. 
for every $L$-formula $\phi (x_1 ,...,x_n)$ and $a_1 ,...,a_n \in N$ 
the values of $\phi (a_1 ,...,a_n )$ in $N$ and in $M$ are the same.}  

\bigskip 

In the case of discrete $\{ 0,1\}$-metric this theorem becomes 
the standard L\"{o}wenheim-Skolem theorem.

\begin{remark} \label{FaHaSh} 
{\em 
Let $L$ be a continuous signature. 
Following Section 4.2 of \cite{fhs} we define a topology 
on $L$-formulas relative to a given continuous theory $T$. 
For $n$-ary formulas $\phi$ and $\psi$ of the same sort set 
$$ 
{\bf d}^T_{\bar{x}} (\phi ,\psi) = \mathsf{sup} \{ |\phi (\bar{a}) -\psi (\bar{a} )|: \bar{a} \in M, M\models T\} . 
$$ 
The function ${\bf d}^T_{\bar{x}}$ is a pseudometric. 
The language $L$ is called {\em separable} 
with respect to $T$ if for any tuple $\bar{x}$ 
the density character of ${\bf d}^T_{\bar{x}}$ is countable. 
By Proposition 4.5 of \cite{fhs} in this case for every 
$L$-model $M\models T$ the set of all interpretations of 
$L$-formulas in $M$ is separable in the uniform topology. 
By Theorem 4.6 of \cite{fhs} if in the formulation of 
the L\"{o}wenheim-Skolem theorem we replace the assumption 
$|L|\le \kappa$ by the condition that $L$ is separable 
with respect to $Th^c (M)$ then the statement 
the L\"{o}wenheim-Skolem theorem also holds for $\kappa =\aleph_0$.   

As we have already noticed the $\mathbb{Q}[i]$-subsignature 
of the language of Hilbert spaces is dense in the standard signature.  
Thus the original language 
of Hilbert spaces is  separable with respect to the theory of Hilbert spaces. 
In particular we may apply the approach of \cite{fhs} in this important case.}  
\end{remark} 

\bigskip

The following corollary of the L\"{o}wenheim-Skolem theorem 
is obvious. 
It in particular states that if a class ${\cal K}$ 
of first-order/continuous-metric 
structures is logically analyzable, then ${\cal K}$ is bountiful.   

\begin{corollary} \label{LS-cor} 
Let $\mathcal{K}$ be a class 
of first-order/continuous-metric structures. 
Assume that there is a family of classes 
$\mathcal{K}_{\alpha}$, $\alpha \in I$, 
in expanded languages with possibly new sorts, so that 
\begin{itemize}   
\item every $G\in \mathcal{K}$ has an expansion $\hat{G}$ 
belonging to one of these classes $\mathcal{K}_{\alpha}$,   
\item the reduct of any elementary substructure of 
such $\hat{G}$ belongs to $\mathcal{K}$.   
\end{itemize} 
Then $\mathcal{K}$ is bountiful. 
\end{corollary}

The following proposition was already mentioned in 
the introduction as a kind of motivation of this paper. 
It can be also considered as a demonstration of 
our method in the easiest form. 
The part of the proof concerning 
proerty {\bf F}$\mathbb{R}$ can be slightly unclear 
for an inexperienced reader. 
Some helpful preliminaries  can be
found in the beginning of Section 3.

\begin{proposition} \label{Bounti} 
(a) Classes of discrete members of {\bf OB}, {\bf FH} and ${\bf F}\mathbb{R}$ 
are not bountiful. 

(b) Classes of discrete members of non-{\bf OB}, non-{\bf FH} 
and non-${\bf F}\mathbb{R}$ are bountiful. 

Moreover statement (a) and (b) also hold in the case of metric groups 
which are continuous structures.  
\end{proposition}  

{\em Proof.} (a) This follows from the theorem of 
Y. de Cornulier that for any finite perfect group 
$F$ and an infinite $I$ the power $F^I$ is strongly bounded 
and the fact that any countable subgroup of   
$F^I$ has cofinality $\omega$, i.e. does not belong to  
${\bf FH} \cup {\bf F}\mathbb{R}$, see \cite{BHV}, \cite{dC} and \cite{pestov16}. 
In particular such a subgroup has an isometric action on a metric 
space with unbounded orbits. 
 
(b) Let us consider the case of groups which are 
first-order structures of the class non-{\bf OB}. 
Let $C\subseteq  H$ with $H\in$ non-{\bf OB}. 
Take an action of $H$ on a metric space $X$ with an unbounded orbit.  
Extend $C$ by countably many elements witnessing this unboundedness. 
By the L\"{o}wenheim-Skolem  theorem there is a subgroup $K$ of $H$ 
which is an elementary substructure of $H$ in the expanded 
(by constants) language such that $C\subseteq  K$ and 
$|C|=|K|$. 
Since this subgroup contains the distinguished constants 
its action on $X$ has unbounded orbits. 

The case of non-{\bf FH} is identical. 
Indeed, by Theorem 2.2.9 of \cite{BHV} property non-{\bf FH} 
is equivalent to existence of an affine isometric action of $G$ 
on a real  Hilbert space so that all (equivalently some) orbits are unbounded.  

In the case of non-${\bf F}\mathbb{R}$ one should use the resut of 
\cite{tits} that each isometry of a real tree is elliptic (i.e. fixes a point) 
or hyperbolic (i.e. acts as a non-trivial translation of some line) and 
every action of a group of elliptic isometries on a tree 
has a fixed point or a fixed end. 
We remind the reader that an end is an equivalence class of half-lines 
under the equivalence relation of having a common half-line. 
The latter statement is a straightforward generalization of 
Exercise 2 of Section 6.5 from \cite{serre} (see Section 1 of \cite{CulMor} 
for the case of $\mathbb{R}$-trees).  
When $H$ has a hyberbolic isometry $h$ all orbits of $\langle h\rangle$ 
are unbounded. 
Thus we apply the argument above arranging  that $K$ has a hyberbolic isometry too. 
When $H$ consists of elliptic isometries and does not have a global fixed point, 
then find a half-line, say $L$, representing the fixed end and take a cofinal sequence 
$c_1 , c_2 , \ldots $ on it. 
For each $c_i$ choose $h_i \in H$ which does not fix $c_1 , \ldots ,c_i$. 
Then  any subgroup of $H$ containing all $h_i$ 
does not have a global fixed point. 
Indeed, if $c$ is such a point, then 
consider segments $[c,c_n ]$. 
All of them are decomposed into a segment belonging to $L$ 
and the unique bridge from $c$ to $L$. 
Thus there is $n_0$ such that for all $n>n_0$ segments 
$[c, c_n ]$ contain $c_{n_0}$.  
This contradicts the assumption that 
$c_{n_0}$ is not fixed by $h_n$ with $n>n_0$.  
 
This argument also works in the case of continuous logic. 
Since first-order structures can be viewed as continuous ones 
there is no need for a continuous version of (a). 
We only mention that when one wants to have 
a non-discrete examples witnessing (a) one can have this by adding  
a compact group as a direct summand. 
$\Box$


\section{Non-amenability vs negation of (T)} 

\paragraph{ A.  Introduction.} 
It is well-known that closed subgroups of amenable 
locally compact groups are amenable. 
This in particular implies that the class of 
amenable locally compact groups which are continuous 
metric structures is bountiful. 
Indeed applying the continuous L\"{o}wenheim-Skolem theorem 
we see that  for any infinite metric group $H$ with a subset 
$C\subseteq H$ where $H$ is amenable and locally compact  
there is an elementary submodel $K$ of $H$ 
such that $C\subseteq  K$ and the density character of $C$  
coincides with the density character of $K$. 
Since $K$ is closed in $H$ it is amenable too. 

\begin{remark} 
{\em It is worth noting that the class 
of all discrete amenable groups is not axiomatizable: 
there are locally finite countable groups having 
elementary extensions containing free groups. }
\end{remark}

Since non-compact amenable locally compact groups  
do not satisfy property {\bf (T)} of Kazhdan 
the argument above suggests that the class 
of non-{\bf (T)} locally compact groups is bountiful too. 
Corollary \ref{bounti-T} below 
is a confirmation of a uniform version of 
this suggestion.   
The main result of this section Theoren \ref{nT} 
shows that in the context of continuous logic 
the class of locally compact groups with 
property non-{\bf (T)} is logically analyzable. 

We apply methods anounced in the introduction. 
The case of non-{\bf (T)} is relatively easy, because 
we only need to consider group actions on Hilbert spaces 
which preserve $n$-balls of $0$. 
It can be considered as a warm up before more 
difficult cases in Section 3. 
In the rest of part A of this section 
we give necessary algebraic definitions. 

Let a topological group $G$ have a continuous unitary 
representation on a complex Hilbert space ${\bf H}$.  
A closed subset $Q\subset G$ 
has an {\bf almost $\varepsilon$-invariant unit vector} in ${\bf H}$ if 
$$ 
\mbox{ there exists }v\in {\bf H} \mbox{ such that } 
\mathsf{sup}_{x\in Q} \parallel x\circ v - v\parallel < \varepsilon
\mbox{ and } \parallel v\parallel =1.  
$$ 
A closed subset $Q$ of the group $G$ is called a {\bf Kazhdan set} 
if there is $\varepsilon >0$ with the following property: 
for every unitary representation of $G$ on a Hilbert space 
where $Q$ has an almost $\varepsilon$-invariant unit vector  
there is a non-zero $G$-invariant vector.  
If the group $G$ has a compact Kahdan subset then 
it is said that $G$ {\bf has property ${\bf (T)}$ of Kazhdan}. 

Proposition 1.2.1 of \cite{BHV} states that the group $G$ 
has  property  ${\bf (T)}$ of Kazhdan if and only 
if any unitary representation of $G$ which weakly contains 
the unit representation of $G$ in $\mathbb{C}$ 
has a fixed unit vector. 

By Corollary F.1.5 of \cite{BHV} the property that 
the unit representation of $G$ in $\mathbb{C}$ 
is {\bf almost contained} in a unitary representation $\pi$ 
of $G$ (this is denoted by $1_{G} \prec \pi$) is 
equivalent to the property that for every compact subset 
$Q$ of $G$ and every $\varepsilon >0$  
the set $Q$ has an almost $\varepsilon$-invariant 
unit vector with respect to $\pi$. 

The following example shows that in the first-order logic 
property ${\bf (T)}$ is not elementary. 

\begin{example} 
{\em Let $n>2$. 
According Example 1.7.4 of \cite{BHV} the group 
$SL_n (\mathbb{Z})$ has property ${\bf (T)}$. 
Let $G$ be a countable elementary extension of 
$SL_n (\mathbb{Z})$ which is not finitely generated. 
Then by Theorem 1.3.1 of \cite{BHV} the group 
$G$ does not have ${\bf (T)}$. 
} 
\end{example} 

\paragraph{B.  Unitary representations in continuous logic.}
In order to treat analyzability question  in the class of 
locally compact groups satisfying some uniform version 
of property non-${\bf (T)}$ we need the preliminaries 
of continuous model theory of Hilbert spaces from Section 1.1. 
Moreover since we want to consider unitary representations of 
metric groups $G$ in continuous logic we should 
fix continuity moduli for the corresponding binary functions 
$G\times B_n \rightarrow B_n$ induced by $G$-actions on 
metric balls of the corresponding Hilbert space. 

This is why we have to consider uniformly  
continuous versions of the notion of a Kazhdan set. 
We define it as follows. 

\begin{definition} 
Let $G$ be a metric group of diameter $\le 1$ 
which is a continuous structure in the language  
$(d,\cdot ,^{-1},1)$.  
Let $\mathcal{F} = \{ F_1 , F_2 , ... \}$ be a family 
of continuity moduli for the $G$-variables of 
continuous function $G \times B_i \rightarrow B_i$. 

We call a closed subset $Q$ of the group $G$ 
an $\mathcal{F} $-Kazhdan set if there is $\varepsilon$ 
with the following property: 
every $\mathcal{F} $-continuous unitary representation of $G$ 
on a Hilbert space with almost $(Q,\varepsilon )$-invariant unit vectors 
also has a non-zero invariant vector.  
\end{definition} 

It is clear that for any family of continuity 
moduli $\mathcal{F}$ a subset $Q\subset G$ 
is $\mathcal{F}$-Kazhdan if it is Kazhdan. 
We will say that $G$ has property 
$\mathcal{F}$-{\bf non-(T)} if $G$ does not have 
a compact $\mathcal{F}$-Kazhdan subset.

To study such actions in continuous logic  
let us consider a class of many-sorted continuous 
metric structures which consist of groups $G$ 
together with metric structures 
of complex Hilbert spaces  
$$ 
(d, \cdot , ^{-1}, 1 ) \cup 
(\{ B_n\}_{n\in \omega} ,0,\{ I_{mn} \}_{m<n} ,
\{ \lambda_c \}_{c\in\mathbb{C}}, +,-,\langle \rangle_{Re} , \langle \rangle_{Im} ).
$$
Such a structure  also contains 
a binary operation $\circ$ of an action which is 
defined by a family of appropriate maps  
$G \times B_m \rightarrow B_{m}$ 
(in fact $\circ$ is presented by a sequence of functions 
$\circ_m$ which agree with respect to all $I_{mn}$). 
When we add the obvious continuous $\mathsf{sup}$-axioms that 
the action is linear and unitary, we obtain an axiomatizable 
class $\mathcal{K}_{GH}$. 
Given unitary action of $G$ on ${\bf H}$
we denote by $A(G,{\bf H})$  the member of $\mathcal{K}_{GH}$ 
which is obtained from this action.  
When we fix continuity moduli, say 
$\mathcal{F} =\{ F_1 ,F_2 ,\ldots \}$, for the $G$-variables 
of the operations $G \times B_m \rightarrow B_{m}$ 
 we denote by 
$\mathcal{K}_{GH}(\mathcal{F} )$ the 
corresponding subclass of $\mathcal{K}_{GH}$. 

\begin{definition} \label{K-delta}  
{\bf  The class}  
$\bigcup \{ \mathcal{K}_{\delta}(\mathcal{F}): \delta \in (0,1) \cap \mathbb{Q}\}$.  
Let 
$\mathcal{K}_{\delta}(\mathcal{F} )$ 
be the subclass of $\mathcal{K}_{GH}(\mathcal{F} )$ 
axiomatizable by all axioms of the following form
$$ 
\mathsf{sup}_{x_1 ,...,x_k \in G} \mathsf{inf}_{v\in B_{m}} 
\mathsf{sup}_{x\in \bigcup x_{i}K_{\delta}} 
\mathsf{max}(\parallel x\circ v - v\parallel \dot- \mbox{ } {1\over n}, 
|1- \parallel v\parallel \mbox{ } |)=0 ,
$$
$$
\mbox{ where } k,m,n\in\omega \setminus \{ 0\} 
\mbox{ and } K_{\delta} = \{ g\in G: d(1,g)\le \delta \}. 
$$ 
\end{definition} 

It is easy to see that the axiom of 
Definition \ref{K-delta} implies that 
each finite union $\bigcup^{k}_{i = 1} g_i K_{\delta}$ 
has an almost ${1\over n}$-invariant unit 
vector in ${\bf H}$. 
To see that it can be written 
by a formula of continuous logic note that  
$\mathsf{sup}_{x\in \bigcup x_{i}K_{\delta}}$
can be replaced by $\mathsf{sup}_{x}$ 
with simultaneuos including of the quantifier-free 
part together with  $\mathsf{max}(\delta\dot{-} d(x,x_i ) : 1\le i\le k)$
into the corresponding $\mathsf{min}$-formula. 

In fact the following theorem shows that the class  of 
$\mathcal{F}$-non-{\bf (T)} locally compact groups 
satisfies the conditions of Corollary \ref{LS-cor}. 

\begin{theorem} \label{nT} 
Let $\mathcal{F} = \{ F_1 , F_2 , ... \}$ be 
a family of continuity moduli for $G$-variables 
of continuous function $G \times B_i \rightarrow B_i$. 

(a) In the class of all unitary $\mathcal{F}$-representations 
of locally compact metric groups the condition of almost containing 
the unit representation $1_{G}$ coincides with 
with the condition of having expansions of the form 
$A(G,{\bf H})$ which are members of 
$\bigcup \{ \mathcal{K}_{\delta}(\mathcal{F}): \delta \in (0,1) \cap \mathbb{Q}\}$. 

(b) If for every compact subset $Q$ of a locally compact 
metric group $(G,d)$ and every $\varepsilon >0$ 
there is an expansion of $G$ to a structure from 
$\mathcal{K}_{GH}(\mathcal{F} )$ 
with a $(Q,\varepsilon )$-almost invariant unit vector 
but without non-zero invariant vectors, then 
$G$ has a unitary $\mathcal{F}$-representation  
which almost contains the unit representation 
$1_{G}$ but does not fix any vector of norm $1$. 
Moreover any elementary substructure of the corresponding 
structure $A(G,{\bf H})$ 
is of the form $A(G_0 ,{\bf H}_0 )$, 
where $G_0 \preceq G$, ${\bf H}_0 \preceq {\bf H}$, 
and also almost contains the unit representation 
$1_{G_0}$ but does not fix any vector of norm $1$. 
\end{theorem} 

{\em Proof.} 
(a) Let $G$ be a locally compact metric group 
and let the ball 
$K_{\delta}= \{ g\in G: d(g,1)\le \delta \}\subseteq G$ 
be compact. 
If a unitary $\mathcal{F}$-representation of $G$ 
almost contains the unit representation 
$1_{G}$, then considering it as a structure 
$A(G,{\bf H})$ we see that this structure belongs to 
$\mathcal{K}_{\delta}(\mathcal{F})$. 

On the other hand if some structure of the form 
$A(G,{\bf H})$ belongs to 
$\mathcal{K}_{\varepsilon}(\mathcal{F})$, 
then assuming that $\varepsilon \le\delta$ 
we easily see that the corresponding representation 
almost contains $1_{G}$. 
If $\delta <\varepsilon$, then $K_{\varepsilon}$ 
may be non-compact. 
However since $K_{\delta} \subseteq K_{\varepsilon}$ 
any compact subset of $G$ belongs to a finite union 
of sets of the form $xK_{\varepsilon}$. 
Thus the axioms of $\mathcal{K}_{\varepsilon}(\mathcal{F})$ 
state that the corresponding structure $A(G,{\bf H})$ 
defines a representation almost containing $1_{G}$. 

(b) Choose $\delta >0$ so that 
the $\delta$-ball 
$K = \{ g\in G: d(g,1)\le \delta \}$ 
in $G$ is compact. 
To see that the group $G$ has a required expansion  
in $\mathcal{K}_{\delta}(\mathcal{F} )$  
we apply the following standard argument 
(see Proposition 1.2.1 from \cite{BHV}). 
For every finite union $\bigcup g_i K$ 
and every $n$ fix a unitary representation 
$\pi_{n,\bar{g}}$ of $G$ without 
non-zero invariant vectors and with a unit vector 
which is $\frac{1}{n}$-invariant with respect to 
$\bigcup g_i K$. 
Then the direct sum of these representations 
almost contains $1_G$. 
Indeed, since every compact subset of $G$ is contained in 
some finite union $\bigcup g_i K$ we see that  
for every compact subset $Q\subset G$ and every 
$\varepsilon >0$ the representation 
has a $(Q,\varepsilon )$-almost invariant unit vector. 
It is clear that there are no non-zero $G$-invariant vectors. 
Let us denote by $M$ the corresponding structure 
from $\mathcal{K}_{GH}(\mathcal{F} )$. 

To see the last assertion of part (b) note 
that since the condition $d(g,1)\le \delta$ 
defines a totally bounded complete subset 
in any elementary extension of $G$, the set $K$ 
above is a definable subset of $acl(\emptyset )$.  

Let $M_0 \preceq M$ and $G_0$ be the sort of $M_0$ 
corresponding to $G$. 
Since the existence of an invariant unit vector 
can be written by a continuous formula we see that 
$M_0$ does not have such a vector. 

It remains to verify that for any compact subset 
$D\subset G_0$ and any $\varepsilon >0$ 
the representation $M_0$ always has 
a $(D,\varepsilon )$-almost invariant unit vector. 
To see this note that since $G_0 \prec G$ and $K$ is 
compact and algebraic, 
the ball $\{ g\in G_0: d(g,1)\le \delta \}\subset G_0$ 
is a compact neighbourhood of 1 which coincides with $K$. 
In particular $D$ is contained in a finite union of 
sets of the form $gK$. 
The rest follows  from the conditions that 
$M_0 \in {\cal K}_{\delta}({\cal F})$ and $G_0 \prec G$. 
$\Box$ 

\bigskip 

\paragraph{C. Comments.} 
In Theorem \ref{nT} we cannot 
axiomatize the class of unitary 
$\mathcal{F}$-representations  
$A(G,{\bf H})$ without 
fixed unit vectors 
(it cannot be done in 
continuous logic). 
Thus the definition of logical analyzability 
is satisfied for property 
$\mathcal{F}$-non-{\bf (T)} 
in a slightly weaker form. 

On the other hand applying Corollary \ref{LS-cor}  
we obtain the following corollary. 

\begin{corollary} \label{bounti-T} 
Locally compact 
metric groups which  have property 
$\mathcal{F}$-non-{\bf (T)}  
form a bountiful class. 
\end{corollary}

 It is an open question if 
the statement of Corollary \ref{bounti-T} 
holds for metric groups which are not 
locally compact.

\begin{remark} \label{QT} 
{\em The author thinks that the following question 
is basic in this topic:}  
\begin{quote}  
Is property {\bf (T)} bountiful in the class of all metric groups?   
\end{quote}  
{\em Analyzing typical examples of groups with  
Kazhdan property {\bf (T)} (for example in \cite{pestov16}) 
it seems likely that bountifulness of {\bf (T)}  
is connected with the following question: }
\begin{quote}
Does an elementary substructure of a discrete group 
with property {\bf (T)} also have  
property {\bf (T)}? 
\end{quote} 
{\em It is natural to consider this question in 
the case of linear groups, where property {\bf (T)} 
and elementary equivalence are actively studied, see 
\cite{EJZ} and \cite{Bunina}. }  
\end{remark}

\begin{remark} 
{\em Since in the locally compact case 
non-compact groups with property {\bf (T)} 
are not amenable the following question 
seems related to Remark \ref{QT}: } 
\begin{quote}  
Is the class of all non-amenable metric groups bountiful?  
\end{quote} 
\end{remark} 

One of definitions of non-amenability says that 
a topological group is  non-amenable if 
there is a locally convex topological 
vector space $V$ and a continuous affine 
representation of $G$ on $V$ such that  
some non-empty invariant convex compact subset $K$ of $V$ 
does not contain a $G$-fixed point (\cite{BHV}, Theorem G.1.7).  

This statement cannot be expressed in logic 
because the notion of locally convex topological 
vector spaces is not logically formalizable. 
On the other hand it is easy to see that if 
we restrict ourself just by linear representations 
on normed/metric vector spaces we obtain a bountiful 
property which is stronger than non-amenability. 
We call this property {\bf strong non-FP}. 

Considering in continuous logic 
linear $G$-representations on  
metric vector spaces $(V,d)$  
we fix continuity moduli for the corresponding binary functions 
$G\times B_n \rightarrow B_n$ induced by the action on 
metric balls $B_n =\{ v\in V: d(0,v) \le n\}$. 
Since the action is not necessary isometric, we now 
need continuity moduli for $B_n$-variables too. 
Thus we define an uniform version strong {\bf non-FP} as follows. 

\begin{definition} 
Let $G$ be a metric group of a bounded diameter 
which is a continuous structure in the language  
$(d,\cdot ,^{-1},1)$.  
Let $\mathcal{F} = \{ F_1 , F_2 , ... \}$ be a family of continuity 
moduli for continuous functions $G \times B_i \rightarrow B_i$. 

We say that $G$ has the strong $\mathcal{F}$-non-{\bf FP} 
property if there is a metric  
vector space $V$ and an $\mathcal{F} $-continuous 
linear representation of $G$ on $V$ such that  
some non-empty invariant convex compact subset $K$ of $V$ 
does not have $G$-fixed points. 
\end{definition}

To realize our approach in this case we should consider 
the class of many-sorted continuous 
metric structures which consists of groups $G$ 
together with metric structures 
of metric vector spaces  
$$ 
(G, d, \cdot , ^{-1}, 1 ) \cup 
(\{ B_n\}_{n\in \omega} ,\{ I_{mn} \}_{m<n} ,
d, 0, \{ \lambda_r \}_{r\in\mathbb{R}}, +,- ) \cup (K, d,I) 
$$ 
and the new sort $K$ corresponding 
to convex compact subspace of $V$. 
It is mapped by $I$ into $B_1$ so that 
$I$ preserves the metric.    

We use the property that when $K$ is compact 
for every natural $n$ there is a number $k_n$ 
such that any subset of $K$ of size $k_n$ 
contains a pair of distance $<\frac{1}{n}$. 
Express this property by a continuous formula, 
say $\phi_n$. 
Note that any  many-sorted 
structure as above which satisfies some 
family of the form 
$\{ \phi_n : n\in \omega \setminus \{ 0\}\}$,  
has the sort $K$ algebraic from the 
point of view of continuous logic, \cite{BYBHU}. 
This implies that any elementary substructure 
has $K$ as a compact sort. 
Now it is easy to verify that the strategy of 
Theorem \ref{nT} and Corollary \ref{bounti-T}
works in this context too. 
In particular we have:

\begin{quote} 
{\em The class of metric groups with 
the strong $\mathcal{F}$-non-{\bf FP} property is 
bountiful.} 
\end{quote}

\begin{remark} 
{\em Let $G$ be a metric group of diameter 1  
which is a continuous structure in the language  
$(d,\cdot ,^{-1},1)$.  
Let $F$ be a continuity modulus for continuous 
functions $G \times B \rightarrow B$ where $B$ is 
a metric space of diameter 1. 

We say that $G$ is $F$-{\bf non-extremely amenable} 
if there is a compact metric space $B$ of diameter 1 
and an $F$-continuous action of $G$ on $B$ which 
does not have a $G$-fixed point in $B$. 
Using the arguments above 
it is easy to see that the following statement holds.}\\ 
The class of metric groups which are 
$F$-non-extremely amenable is bountiful. 
\end{remark}


\section{Unbounded actions} 

In this section we consider actions which 
do not preserve $n$-balls 
of metric spaces. 
This situation is more complicated 
than the one of Section 2.

The following material is standard, \cite{EsKi}. 
Let $(X,d)$ be a metric space.  
It is called {\bf pointed} if we fix a point from $X$. 
A {\bf geodesic path} joining $x\in X$ to $y\in X$ 
is an isometric map $\rho$ from some closed 
interval $[0,l]\subset \mathbb{R}$ to $X$ 
such that $\rho(0) =x$ and $\rho (l) = y$. 
Let $[x,y]=\rho ([0,l])$. 
The space $(X,d)$ is called {\bf uniquely geodesic} 
is there is exactly one geodesic path joining $x$ and $y$ for each 
$x,y\in X$. 
Note that the Hilbert space ${\bf H}$ is uniquely geodesic. 
A uniquiely geodesic space is called an $\mathbb{R}$-{\bf tree} 
if for any $x,y,z$ the condition $[y,x]\cap [x,z] = \{ x\}$ implies 
$[y,x] \cup [x,z] = [y,z]$. 
A subset $S\subseteq X$ is {\bf convex} if 
$(\forall x,y \in S) [x,y] \subseteq S$. 

Assume that $X$ is an $\mathbb{R}$-tree. 
Then a convex subset is also called a {\bf subtree}. 
Given $x,y,z \in X$, 
there is a unique element  $c \in [x,y]\cap [y,z] \cap [z,x]$,
called the {\bf median} of $x,y,z$. 
When $c\notin\{x,y,z\}$, the subtree $[x,y]\cup [x,z]\cup [y,z]$ 
is called a {\bf tripod}.
A {\bf line} is a convex subset containing no tripod and maximal for inclusion. 
The {\bf betweenness} relation $B$ of $X$ is the ternary relation
$B(x;y,z)$ defined by $x\in (y,z)$.
For any line $L$ and a point $a\in L$ the relation $\neg B(a; y,z)$ 
is an equivalence relation on $L\setminus \{ a\}$. 
An equivalence class of this relation together with $a$ 
is called a {\bf half-line}. 

We say that half-lines $L_1$ and $L_2$ are equivalent if 
$L_1 \cap L_2$ contains a half-line. 
An {\bf end} is an equivalence class of this relation.

\bigskip 

We consider complete uniquely geodesic metric spaces as a 
many-sorted metric structures of $n$-balls of pointed spaces  
$$ 
(\{ B_n\}_{n\in \omega} ,0,\{ I_{mn} \}_{m<n} ,d ).
$$
where $0$ denotes the fixed point and other symbols 
are interpreted in the natural way. 

It is shown in \cite{carlisle} 
that the class of pointed real trees is axiomatizable 
in continuous logic by axioms of $0$-hyperbolicity and 
the approximate midpoint property. 
We remind the reader that a metric space 
$(X,d)$ has the {\bf approximate midpoint property}   
if for any $x,y\in X$ and any rational $\varepsilon >0$ 
there exists $z\in X$ such that  
$$
|d(x,z) - d(x,y)/2 | \le \varepsilon \mbox{ and } 
|d(y,z) - d(x,y)/2 | \le \varepsilon .   
$$ 
The space $(X,d)$ is called 0-{\bf hyperbolic} if for any 
$x,y,z,w \in X$ and any rational $\varepsilon >0$  
$$ 
\mathsf{min} ( (x\cdot y)_w , (y\cdot z)_w ) \le (x\cdot z)_w +\varepsilon 
\mbox{ , where } 
$$  
$$ 
(x\cdot y)_w  =\frac{1}{2}( d(x ,w) + d(w,y) - d(x,y) ) \mbox{ (Gromov product). } 
$$

\bigskip 

We now give necessary 
information concerning  isometric actions
of groups on real trees. 
This is taken from \cite{CulMor}, 
\cite{serre} and \cite{tits}. 

Let $T$ be a real tree. 
Assume that $G$ has an isometric action on an $\mathbb{R}$-tree $T$.
We say that $g\in G$ is {\bf elliptic} if it has a fixed point, 
and {\bf hyperbolic} otherwise.

\begin{lemma} \label{lem_nnest}
Let $G$ be a group with an isometric action on an $\mathbb{R}$-tree $T$.

  \begin{itemize}
  \item If $g\in G$ is elliptic, its set of fix points $T^g$ is a closed convex subset.
  \item\label{nnest_it2} If $g$ is hyperbolic, there exists a unique line $L_g$ preserved by $g$; moreover,
$g$ acts on $L_{g}$ by a translation. 
  \item If $g$ is hyperbolic, then for any $p \in T$, $[p,g(p)]$ meets $L_{g}$ 
and $[p,g(p)] \cap  L_{g} = [q,g(q)]$ for some $q \in L_{g}$. 
  \item If $g$ and $h\in G$ are elliptic and $T^g \cap T^h =\emptyset$, then 
  $gh$ is hyperbolic. 
  \end{itemize} 
\end{lemma}

When $g$ is hyperbolic, $L_g$ is called the {\bf axis} of $g$.

\subsection{Actions and continuity moduli} 

\paragraph{A. Assumptions.}
Let a metric group   
$(G, \cdot , ^{-1}, d )$ act on $(X,d)$ by isometries, 
where $(X,d)$ is a uniquely geodesic space.  
When we consider this situation we always take 
the following assumptions. 
They are necessary to present the situation in continuous logic. 

We firstly fix a point $0\in X$, present $X$ as 
the union of $n$-balls of $0$, and 
assume the existence of a function 
$\mathsf{gth}: \mathbb{N} \times [0,1) \rightarrow \mathbb{N}$ 
such that for every natural $m$ 
if $g\in G$  satisfies $d(1,g)=\delta < 1$, then 
$g$ takes the ball $B_m$ into the ball 
$B_{\mathsf{gth}(n,\delta )}$. 
We will assume that $\mathsf{gth}$ is 
increasing with respect to the first argument. 
Note that when $(G,d)$ is discrete with 
with the $\{ 0,1\}$-metric we can take 
$\mathsf{gth}(n,\delta ) =\mathsf{id}(n)$. 

Note that for any pair $m<n$ the action defines a partial map 
$G\times B_m \rightarrow B_n$ which is uniformly continuous 
with respect to the $B_m$-argument.    
In order to present this action in continuous logic we assume that 
this function is uniformly continuous 
with respect to the $G$-argument too. 
To do this in a formal way 
let us define ternary predicates 
$$
\circ_{mn}(g,x,y):G\times B_m \times B_m \rightarrow [0, m+n ]  
\mbox{ , } m\le n\in \mathbb{N} ,    
$$ 
as follows. 
\begin{definition}  \label{circ} 
$$
\circ_{mn}(g,x,y) = \mathsf{length} ([g\circ x ,y]\cap B_n ), 
\mbox{ where }x\in B_m \mbox{ and } y \in B_m . 
$$ 
\end{definition} 
In this formula the length is defined with respect 
to the metric of $X$. 
Note that for $x,y\in B_m$ with $g\circ x \in B_n$ 
we have $\circ_{mn}(g,x,y) = d(g\circ x ,y)$. 
We now formulate our basic assumption. 

\begin{definition} \label{def-cont} 
Let $\mathsf{ort} : \mathbb{N} \rightarrow \mathbb{N}$ 
be a function such that for any $m,n \in \mathbb{N}$ 
with $\mathsf{ort}(m) <n$ there is a continuity 
modulus for the $G$-variable of the 
predicate $\circ_{mn}$. 

We fix such a modulus and denote it by $\gamma^{m,n}_G$. 
\end{definition}

\begin{remark} 
{\em It is worth mentioning here that the assumption 
of uniform continuity of the unrestricted function 
$\circ :G\times X \rightarrow X$ strongly 
trivialize the situation. 
Indeed, assume that $X$ is an $\mathbb{R}$-tree and 
$g\in G$ acts on $X$ with non-fixed points. 
If $g$ is ellipic, fixes $x_0$  and takes some $x$ to $y\not = x$, 
then a half-line starting at $x_0$ and containing $x$ is 
taken by $g$ to a half-line containing $y$. 
If the first half-line is unbounded the set of 
values $d(z, g\circ z )$ for elements $z$ 
of this half-line is unbounded too. 
A similar argument works in the case of hyperbolic $g$. 
In this case one should consider 
half-lines not cofinal with the axis of $g$.  
We see that for non-discrete $(G,d)$ no condition 
of the form $d(1,g) < \delta$ bounds 
the set of values $d(z, g\circ z )$ if 
$X$ has infinitely many pairwise nonequivalent 
unbounded half-lines. 
In particular the action is not uniformly 
continuous as a binary function.}
\end{remark} 

We now discuss the assumptions of existence of  
functions $\mathsf{gth}$ and $\mathsf{ort}$ 
in the case of $\mathbb{R}$-trees and in the case 
of real ${\bf H}$ separately. 
We will see below that in the former case 
some natural geometric condition guarantees that 
the identity function $\mathsf{id}$ works for  
$\mathsf{ort}$. 
In the latter case function $\mathsf{id}$ is not relevant. 
However we will give a completely satisfactory answer in this case: 
{\em the only assumption that all partial maps 
$G\times B_m \rightarrow B_n$ defined by the action
on ${\bf H}$ are uniformly continuous implies 
the existence of a function $\mathsf{ort}$ 
satisfying Definition \ref{def-cont}. }

\bigskip 

\paragraph{B. The case of $\mathbb{R}$-trees.} 
Let a metric group   
$(G, \cdot , ^{-1}, d )$ act on $(X,d)$ by isometries, 
where $(X,d)$ is a pointed $\mathbb{R}$-tree.  
For any pair $m<n$ let $\tilde{\gamma}^{m,n}_{G}$ 
be a continuity modulus of the $G$-variable  
of the partial function $G\times B_m \rightarrow B_n$.  

In Definition \ref{def-hyp} we formulate a condition 
which guarantees the assumption of 
Definition \ref{def-cont} in the case $\mathsf{ort} = \mathsf{id}$.  
This will be proved in Lemma \ref{R-cont}. 
In Remark \ref{R-cont-rem} we will see that 
these conditions are equivalent.

\begin{definition} \label{def-hyp}
We say that the action of $G$ on $X$ has  {\bf hyperbolic continuity moduli} 
if for any $h\in G$, any pair $m<n$ and any rational $q >0$ 
if  $x\in B_m$ and $y\in h(B_m )\setminus B_n$ then   
$$
\forall g \in G 
(d(1,g)\le \tilde{\gamma}^{m,n}_G (q ) \rightarrow |\mathsf{length} ([x,y]\cap B_n ) 
- \mathsf{length} ([x, g\circ y]\cap B_n )| \le q ).   
$$ 
\end{definition} 

In fact this property links the metric of $G$ with the metric of $X$. 
For an illustration consider the case when $X$ is 
a countable simplicial tree of finite valency and $G\le Aut (X)$ 
is under a metric defining the pointwise convergence topology.  
Then the property of hyperbolic continuity moduli obviously holds   
after some correction of $\tilde{\gamma}^{m,n}_G$ if necessary:  
we only arrange that  
for any $g\in G$ with $d(1,g)\le \tilde{\gamma}^{m,n}_G (q)$ 
the set $B_n$ is contained in the set of 
fixed points of $g$. 
In this case we even have 
$$
\mathsf{length} ([x,y]\cap B_n ) = \mathsf{length} ([x, g\circ y]\cap B_n ) .
$$  

\begin{remark} 
{\em  
The property of hyperbolic moduli follows from the condition 
that for any $n$ and any $q>0$ there exists $\delta>0$  
such that 
$$ 
\forall g \in G \forall y\in X 
(d(1,g)\le \delta \wedge y\notin B_n \rightarrow \mathsf{length} ([y, g\circ y]\cap B_n ) \le q ).   
$$ 
Indeed correcting $\tilde{\gamma}^{m,n}_G$ if necessary we 
can replace $\delta$ in this formula by $\tilde{\gamma}^{m,n}_G (q)$. 
Now to verify the condition of Definition \ref{def-hyp} 
take the median, say $c$, of the triple $x,y,g\circ y$. 
If $c\not\in B_n$ then 
$$
\mathsf{length} ([x,y]\cap B_n ) = \mathsf{length} ([x, g\circ y]\cap B_n ) .
$$ 
If $c\in B_n$ then 
$$
|\mathsf{length}([x,y] \cap B_n) - \mathsf{length}([x,g\circ y]\cap B_n)| \le 
$$ 
$$
|\mathsf{length}([y,c] \cap B_n) + \mathsf{length}([c,g\circ y]\cap B_n)| , 
$$ 
and the latter value is not greater than $q$.   
} 
\end{remark}

\begin{lemma} \label{R-cont} 
Assume that a group $G$ acts on $X =\bigcup B_n$ 
by isometries and the action $\circ$ 
has hyperbolic continuity moduli $\tilde{\gamma}^{m,n}_G$ on $G$. 
Let $A(G,X)$ be the corresponding many-sorted 
metric structures of $n$-balls of $0$    
$$ 
(G, \cdot, d) \cup (\{ B_n\}_{n\in \omega} ,0,\{ I_{mn} \}_{m<n} ,d ).
$$
extended by the predicates $\circ_{mn}$, $m<n$, as above. 

Then the predicates $\circ_{mn}$ have continuity moduli 
$\gamma^{m,n}_G = \tilde{\gamma}^{m,n}_G$ for the $G$-variable 
and $\mathsf{id}$ for variables of $X$.  
\end{lemma}

{\em Proof.} 
Let us consider the first $B_m$-variable. 
Assume $d(x,x')<\varepsilon$. 
Then $d(g\circ x ,g\circ x')<\varepsilon$.
If $g\circ x,g\circ x'\in B_n$ then obviously 
$|d(g\circ x,y)- d(g\circ x',y)|< \varepsilon$. 
Let $g\circ x'\not\in B_n$ and let $c$ 
be the median of the triple $y,g\circ x,g\circ x'$. 
Then $[y,g\circ x] = [y,c]\cup [c,g\circ x]$ and 
$[y,g\circ x'] = [y,c]\cup [c,g\circ x']$. 
If $[y,c]$ contains a point of distance $n$ 
from $0$ then $\circ_{mn}(g,x,y) = \circ_{mn}(g,x',y)$. 
If $[y,c]$ does not contain a point 
of distance $n$ from $0$  then 
$$
|\circ_{mn}(g,x,y) -\circ_{mn}(g,x',y)|\le 
|d(g\circ x,c) - d(g\circ x',c)|<\varepsilon. 
$$ 
The case of the second $B_m$-variable is similar. 

In the case of the variable of $G$ assume that 
$d(g,g') \le \tilde{\gamma}^{m,n}_G (\varepsilon )$. 
Then for $g\circ x, g'\circ x \in B_{n}$ we have 
$d(g\circ x, g'\circ x) \le \varepsilon$ and by the triangle inequality, 
$|\circ_{mn} (g,x,y ) - \circ_{mn} (g',x,y )|\le \varepsilon$. 

If $g'\circ x \not\in B_{n}$ then 
$$ 
|\circ_{mn}(g,x,y) -\circ_{mn}(g',x,y)|\le 
$$
$$
|\mathsf{length}([y,g'\circ x] \cap B_n) - \mathsf{length}([y,g(g')^{-1} \circ (g'\circ x)]\cap B_n)| 
$$ 
and the latter value does not exceed 
$\varepsilon$ by hyperbolicity of continuity moduli. 
The remaining case is similar. 
$\Box$ 

\bigskip 

\begin{remark} \label{R-cont-rem} 
{\em Note that Lemma \ref{R-cont} also holds in the opposite direction. }  

Assume that a group $G$ acts on $X =\bigcup B_n$ 
by isometries and the predicates $\circ_{mn}$, $m<n$,
have continuity moduli $\gamma^{m,n}_G$ 
which are also continuity moduli of 
functions $G\times B_m \rightarrow B_n$ 
defined by the action $\circ$.  
Then these continuity moduli satisfy 
the condition of hyperbolicity from   
Definition \ref{def-hyp}. 

{\em In terms of Definition \ref{def-hyp} to see this 
statement it suffices to consider the value}  
$$ 
|\circ_{mn}(h,h^{-1}(y),x) -\circ_{mn}(gh,h^{-1}(y),x)| . 
$$
\end{remark}

Lemma \ref{R-cont} shows that typical $G$-spaces 
which are real trees can be considered as 
continuous metric structures with respect to 
some functions 
$\mathsf{gth}$ and $\mathsf{ort}$. 
Indeed we can define $\mathsf{ort} (m) = m$. 
In Section 3.2 we demonstrate that this approach 
works for logical analyzability of the class  
non-${\bf F}\mathbb{R}$ with respect to actions 
with some $\mathsf{gth}$ and $\mathsf{ort}$ defined.  
In particular 
\begin{itemize} 
\item in Theorem \ref{isoR} we prove that the class of 
structures $A(G,X)$ as in Lemma \ref{R-cont} 
(and with fixed $\mathsf{gth}$ and $\mathsf{ort}$) 
is axiomatizable, 
\item the statement of 
Proposition \ref{Bounti}(b) concerning non-{\bf F}$\mathbb{R}$ 
can be extended to a statement of logical analyzability, 
see Theorem \ref{n-FR}. 
\end{itemize} 
Does this approach work in the case of 
isometric actions on ${\bf H}$? 
It is easy to see that the definition 
$\mathsf{ort}(m) = m$ cannot be 
justified in the case of ${\bf H}$. 
We now discuss this problem and 
show how our approach 
should be modified in this case. 

\bigskip

\paragraph{C. The case of Hilbert spaces.}

Assume that a metric group   
$(G, \cdot , ^{-1}, d )$ acts on the real Hilbert space 
${\bf H}$ by isometries with respect to the metric induced by 
$\parallel \parallel$. 
It is well-known that such isometries are affine transformations 
(\cite{BHV}, Chapter 2).    

We start with the description why in the case 
of Hilbert spaces the ternary predicates $\circ_{mn}$ 
defined above can lose continuity moduli for some $m<n$. 
Let us fix $k,l,m\in \mathbb{N}\setminus \{ 0\}$ 
and $\varepsilon>0$. 
Let $n = m+1$. 
Take a point $y_1$ of norm $m$ and 
a point $y_2 \in B_{n}\setminus B_{n-\frac{\varepsilon}{k}}$.  
Assume that $d(y_1 ,y_2)=l$. 
Let $g\in G$ satisfy $g\circ 0\not\in B_n$,   
$d(g\circ 0,y_2 )= \varepsilon$ and 
$y_2 \in [0,g\circ 0]$. 
Let $g'\circ 0 = y_2$. 
Assuming that $k$, $l$ and $m$ are sufficiently large 
(in particular the $n$-sphere bounding $B_n$ is close to 
be flat in some domain containing $y_1$ and $y_2$) 
we can arrange that 
$|\circ_{mn} (g,0,y_1) - \circ_{mn} (g',0,y_1)|$ 
is close to $1$, while $d(g\circ 0,g'\circ 0) = \varepsilon$.    
In particular we do not have a method of 
defining continuity moduli in this case. 

The following geometric lemma 
shows how the function $\mathsf{ort}$ should be defined.

\begin{lemma} \label{triangle} 
For every natural number $m$ there exists a number 
$n>m$ such that the following statements hold. \\
(a) For any positive $\varepsilon <\frac{1}{2}$ if 
$v,v_1 ,v_2 \in {\bf H}$ satisfy   
$\mathsf{min} (\parallel v_1 \parallel ,\parallel v_2 \parallel ) \geq n$, 
$\parallel v\parallel \le m$, and 
$\parallel v_1 - v_2 \parallel <\varepsilon$, then  
the distance between the points of the $\parallel \parallel$-norm 
$n$, say $v'_1$ and $v'_2$, belonging correspondingly to the segments 
$[v,v_1]$ and $[v,v_2 ]$  
is smaller than $2\varepsilon$. \\ 
(b) 
For any positive $\varepsilon <\frac{1}{2}$ if 
$v,v_1 ,v_2 \in {\bf H}$ satisfy   
$\mathsf{max} (\parallel v_1 \parallel ,\parallel v_2 \parallel ) \le m$, 
$\parallel v\parallel >n$, and 
$\parallel v_1 - v_2 \parallel <\varepsilon$, then  
the distance between the points of the $\parallel \parallel$-norm 
$n$, say $v'_1$ and $v'_2$, belonging correspondingly to the segments 
$[v,v_1]$ and $[v,v_2 ]$  
is smaller than $2\varepsilon$. \\ 
(c) For any positive $\varepsilon <\frac{1}{2}$ if 
$v,v_1 ,v_2 \in {\bf H}$ satisfy   
$\parallel v_1 \parallel \ge n$, $\parallel v_2 \parallel  \le m$, 
$\parallel v\parallel =n$, $v\in [v_1 ,v_2 ]$ and 
$v_1 \in B_{n+\varepsilon}$, then  
$\parallel v_1 - v \parallel \le 2\varepsilon$.
\end{lemma} 

{\em Proof.} 
We prove (a) and (b) simultaneously. 
Let us choose $n$ so that for any $u_1 ,u_2 \in B_m$ 
and any $v_1 ,v_2\not\in B_n$ with 
$d(v_1 ,v_2 )\le \frac{1}{2}$ 
the angle $\alpha$ 
between $\overrightarrow{u_1 v_1 }$ and $\overrightarrow{u_2 v_2 }$ 
is sufficiently small. 
In particular these vectors are close to be collinear 
to $\overrightarrow{0 v_1 }$ and $\overrightarrow{0 v_2 }$. 
Then for $v'_1 ,v'_2$ of the norm $n$ 
so that $v'_1 \in [u_1,v_1 ]$ and $v'_2 \in [u_2,v_2 ]$,
the vector $\overrightarrow{v'_1 v'_2 }$ is close to be  
orthogonal both to $\overrightarrow{u_1 v_1 }$ and $\overrightarrow{u_2 v_2 }$. 
In particular we may assume that when 
$d(u_1 ,u_2 )\le \varepsilon \le \frac{1}{2}$ 
and $d(v_1 ,v_2 )\le \varepsilon \le \frac{1}{2}$ 
we have $d(v'_1 ,v'_2 )\le 2 \varepsilon$.  

Statement (c) follows from the fact that 
$\overrightarrow{v_1 v_2}$ is close to be 
collinear with $\overrightarrow{v_1 0}$. 
$\Box$ 
\bigskip

Let us fix a function 
$\mathsf{ort}: \mathbb{N} \rightarrow \mathbb{N}$ 
which finds for every $m$ 
a number $n$ as in the formulation of Lemma \ref{triangle}. 
Notice that this function only depends on 
geometric properties of ${\bf H}$.  
Let $\mathsf{ort}^H$ be the minimal one.  
We consider ternary predicates 
$$
\circ_{mn}(g,x,y):G\times B_m \times B_m \rightarrow [0, n+m ]  
\mbox{ , } m\in \omega \mbox{ , } \mathsf{ort}^H (m)\le n ,    
$$ 
as in Definition \ref{circ}. 
Let $A(G,{\bf H})$ be the corresponding many-sorted 
metric structures of $n$-balls of $0$    
$$ 
(G, \cdot, d) \cup (\{ B_n\}_{n\in \omega} ,0,\{ I_{mn} \}_{m<n} ,d ).
$$
extended by the predicates $\circ_{mn}$ for $\mathsf{ort}^H (m) <n$.

\begin{lemma} \label{H-cont} 
Let a group $G$ act on ${\bf H} =\bigcup B_n$ 
by affine isometries. 
Assume that for any $m,n$ with $\mathsf{ort}^H (m)<n$ 
there is a continuity modulus $\tilde{\gamma}^{m,n}_G$ 
for the $G$-variable of the partial function 
$G\times B_m \rightarrow B_n$ induced by the action. 
Let $A(G,{\bf H})$ be the corresponding structure 
with all $\circ_{mn}$ where $\mathsf{ort}^H (m)< n$. 
Then the predicates $\circ_{mn}$ have continuity moduli 
$8\tilde{\gamma}^{m,n}_G$ for the $G$-variable and $3\mathsf{id}$ 
for variables of ${\bf H}$.  
\end{lemma}

{\em Proof.} 
The proof is similar to Lemma \ref{R-cont}. 
Let us consider the first $B_m$-variable. 
Assume $d(x,x')<\varepsilon< \frac{1}{2}$. 
Then $d(g\circ x ,g\circ x')<\varepsilon$.
If $g\circ x,g\circ x'\in B_n$ then obviously 
$|d(g\circ x,y)- d(g\circ x',y)|< \varepsilon$. 
Let $g\circ x'\not\in B_n$  
and let $x'' \in [y, g\circ x']$ be the 
point of distance $n$ from $0$. 
Then in the case $g\circ x\in B_n$ Lemma \ref{triangle}(c)  
and the triangle inequality guarantee that $d(g\circ x,x'') \le 3\varepsilon$,  
i.e. $|d(g\circ x,y)- d(x'',y)|< 3\varepsilon$. 
When $g\circ x \not\in B_n$ choose 
$\hat{x} \in [y,g\circ x]$ of distance $n$ from $0$. 
Then by Lemma \ref{triangle}(a) we have 
$|d(\hat{x},y)- d(x'',y)|< 2\varepsilon$ 
which in fact is equivalent to 
$|\circ_{mn}(g,x,y)- \circ_{mn}(g,x',y)|< 2\varepsilon$.

The case of the $B_m$-varable $y$ is similar. 

In the case of the variable of $G$ assume that 
$d(g,g') \le \tilde{\gamma}^{m,n}_G (\varepsilon )$. 
If $g\circ x$, $g'\circ x \in B_n$ then 
$d(g\circ x, g'\circ x) \le \varepsilon$ 
and by the triangle inequality, 
$$
|\circ_{mn} (g,x,y ) - \circ_{mn} (g',x,y )|\le \varepsilon . 
$$ 
Assume that neither $g\circ x$ nor $g'\circ x$ 
belongs to $B_n$. 
The map $g'g^{-1}$ takes $[y, g\circ x ]$ to 
$[g'g^{-1}\circ y, g'\circ x]$. 
Let $y' \in [y, g\circ x ]$ and 
$y'' \in [g'g^{-1}\circ y, g'\circ x]$
be elements of norm $n$. 
Replacing $y'$ by $y''$ if necessary we may assume 
below that $(g'g^{-1})^{-1} \circ y'' \in B_n$. 
Then 
$\parallel y'' -  (g'g^{-1})^{-1} \circ y'' \parallel \le \varepsilon$ 
and 
$\parallel y -  g'g^{-1} \circ y \parallel \le \varepsilon$. 
Since $n - \parallel (g'g^{-1})^{-1} \circ y'' \parallel \le \varepsilon$ 
applying arguments as in Lemma \ref{triangle} 
we easily see that 
$\parallel y' -  (g'g^{-1})^{-1} \circ y'' \parallel \le 2 \varepsilon$ 
and 
$\parallel y' -  y'' \parallel \le 3\varepsilon$. 
In particular the difference between 
$\parallel y - y'\parallel$ and 
$\parallel g'g^{-1}\circ y - y''\parallel$ 
does not exceed $5 \varepsilon$. 
Applying Lemma \ref{triangle} (b) we see that the 
length of the $B_n$-part of the segment $[y, g'\circ x]$ 
does not differ from 
$\parallel g'g^{-1}\circ y - y''\parallel$ 
more than $3 \varepsilon$. 
In particular it does not differ from 
$\parallel y-y'\parallel$ more than $8 \varepsilon$. 
By definition of $\circ_{mn}$ we have 
$$
|\circ_{mn} (g,x,y ) - \circ_{mn} (g',x,y )|\le 8 \varepsilon . 
$$ 
The case when  $|\{ g\circ x , g'\circ x \}\cap B_n |=1$ 
can be arranged by similar arguments. 
$\Box$ 

\bigskip

\subsection{Non-F$\mathbb{R}$-actions}

Let us fix  continuity moduli $\gamma^{m,n}_G$ 
and consider the class of continuous metric structures 
which are unions of continuous groups  
$(G, \cdot , ^{-1}, d )$  
together with many-sorted metric structures  
of $n$-balls of pointed $\mathbb{R}$-trees   
$$ 
{\bf X} = (\{ B_n\}_{n\in \omega} ,0,\{ I_{mn} \}_{m<n} ,d ) 
\mbox{ , where } I_{mn}: B_m\rightarrow B_n ,  
$$ 
and ternary predicates   
$$
\circ_{mn}(g,x,y):G\times B_m \times B_m \rightarrow [0, m+n ]  
\mbox{ , } m\le n\in \omega ,   
$$ 
where $\gamma^{m,n}_G$ is the continuity modulus 
with respect to $G$ and the continuity moduli 
with respect to $B_m$  are equal to $\mathsf{id}$.

We now introduce an axiomatizable class 
$\mathcal{K}_{iso\mathbb{R}}(\mathsf{gth}, \gamma^{m,n}_G)$ 
which contains all structures of the 
form $A(G,X)$ from Lemma \ref{R-cont} 
which correspond to actions with a fixed function 
$\mathsf{gth}: \mathbb{N}\times [0,1) \rightarrow \mathbb{N}$  
(when we assume that $\mathsf{ort} = \mathsf{id}$). 
We will see in Theorem \ref{isoR} that 
the axioms below guarantee that all members of 
$\mathcal{K}_{iso\mathbb{R}}(\mathsf{gth},\gamma^{m,n}_G)$ 
are obtained in this way. 
This is a confirmation that 
our approach works well!

Since our formulas are becoming quite complicated 
we will not mention below the domains of quantifiers 
$\mathsf{inf}$ and $\mathsf{sup}$. 
Usually this is clear from the sorts of the predicates 
appearing in the formula. 

\begin{definition} \label{ClassAx} 
{\em Let $\mathcal{K}_{iso\mathbb{R}}(\mathsf{gth}, \gamma^{m,n}_G)$ 
be the class of structures of the language above which 
satisfy axioms of metric groups with invariant metrics for the 
sort $G$, the axioms of pointed real trees 
and the following three groups of axioms of the relation 
$\circ_{mn}$.   \\ 
{\em 1(a) The value $g\circ x$ is eventually defined inside some $B_n$.} \\ 
$$ 
\mathsf{sup}_g \mathsf{sup}_x \mathsf{min} ( \circ_{ms} (g,x,0)\dot{-} n , 
\delta \dot{-} d(1,g)) = 0,  
$$ 
where $\delta \in [0,1) \cap \mathbb{Q}$  and $\mathsf{gth} (m, \delta )\le n<s$. \\ 
{\em 1(b) Correctness.} 
$$ 
\mathsf{sup}_g \mathsf{sup}_{x,y} \mathsf{min} ( |\circ_{ms} (g,x,y)\dot{-} \circ_{m't} (g, I_{mm'} (x),I_{mm'} (y))| , 
\delta \dot{-} d(1,g)) = 0,  
$$
where $\delta \in [0,1) \cap \mathbb{Q}$, $m<m'$, 
$\mathsf{gth} (m, \delta )<s<t$ and $\mathsf{gth} (m', \delta )<t$. \\ 
{\em 1(c) Approximation of the value $g\circ x$ by points $z$ 
(when $g\circ x \in B_m$). }\\ 
$$
\mathsf{sup}_g \mathsf{sup}_x \mathsf{min} ( m\dot{-}\circ_{mn} (g,x,0 ), \Phi (g,x))= 0 \mbox{ , where }
$$ 
$$
\Phi (g,x) = \mathsf{inf}_z \mathsf{max} (\circ_{mn} (g,x,z ), \mathsf{sup}_u \mathsf{min} (|d(u,z) -\circ_{mn} (g,x,u )|, m\dot{-} d(0,u) )).  
$$ 
{\em 2(a) The predicate $\circ_{mn}$ measures the distance inside $B_m$.  
Triangle inequalities involving $g \circ x$.}   \\ 
$\mathsf{sup}_{g} \mathsf{sup}_{x,y_1,y_2}\mathsf{min} (d(y_1 ,y_2 )\dot{-} (\circ_{mn} (g,x,y_1 )+\circ_{mn} (g,x,y_2 ))$, 

\hspace{8cm} $m\dot{-} \circ_{mn} (g,x,0 ))= 0$. \\   
$\mathsf{sup}_{g} \mathsf{sup}_{x,y_1,y_2}\mathsf{min} (\circ_{mn} (g,x,y_1 )\dot{-} (d(y_1 ,y_2 )+\circ_{mn} (g,x,y_2 ))$,

\hspace{8cm} $m\dot{-} \circ_{mn} (g,x,0 ))= 0$. \\  
{\em 2(b) A version of the triangle inequality when the action $\circ$ is isometric.} \\  
$\mathsf{sup}_{g} \mathsf{sup}_{x_1 ,x_2}\mathsf{min} (\mathsf{sup}_y (d(x_1 ,x_2 ) \dot{-} (\circ_{mn}(g,x_1 ,y)+\circ_{mn}(g,x_2 ,y) )),m\dot{-}\circ_{mn} (g,x_1,0 )$, 

\hspace{8cm} $m\dot{-} \circ_{mn} (g,x_2,0 ))= 0$. \\
{\em 2(c) The action $\circ$ is isometric.}  \\ 
$\mathsf{sup}_{g} \mathsf{sup}_{x_1 ,x_2}\mathsf{min} (\mathsf{inf}_y ( (\circ_{mn}(g,x_1 ,y)+\circ_{mn}(g,x_2 ,y) )\dot{-}d(x_1 ,x_2 )),m\dot{-} \circ_{mn} (g,x_1,0 )$,

\hspace{8cm} $m\dot{-} \circ_{mn} (g,x_2,0 ))= 0$. \\ 
{\em 3(a) The neutral element acts trivially.} \\ 
$\mathsf{sup}_{x ,y}(|\circ_{mn} (1,x,y) - d(x,y)|)= 0$.  \\   
{\em 3(b) The action of $g^{-1}$ inside $B_m$.}  \\
$\mathsf{sup}_{g}\mathsf{sup}_{x ,y}\mathsf{min}( |\circ_{mn} (g,x,y)- \circ_{mn} (g^{-1},y,x)| , 
m\dot{-} \circ_{mn} (g,x,0 ), $

\hspace{8cm} $m\dot{-} \circ_{mn} (g^{-1},y,0 ))= 0$.   \\ 
{\em 3(c) When $g'(x)=z$ and $g(z) = y$, then $gg'(x) = y$.} \\  
$$ 
\mathsf{sup}_{x,y,g,g'}\mathsf{min} (\Phi (g,g',x,y),m\dot{-}\circ_{mn} (g^{-1},y,0 ),m\dot{-}\circ_{mn} (g',x,0 ),m\dot{-}\circ_{mn} (gg',x,0 ))= 0,  
$$ 
$$ 
\mbox{ where } \Phi(g,g',x,y) = \mathsf{sup}_{z} (\circ_{nm} (gg', x,y) \dot{-}(\circ_{nm} (g^{-1} ,y, z) + \circ_{nm} (g' ,x,z))).  
$$ 
}
\end{definition} 

\bigskip

\begin{theorem} \label{isoR} 
(a) Assume that a group $G$ acts on $X =\bigcup B_n$ 
by isometries according to a function $\mathsf{gth}$ 
and the action $\circ$ 
has hyperbolic continuity moduli 
$\{ \gamma^{m,n}_G : m<n \in \omega\}$ on $G$. 
Let $A(G,X)$ be the corresponding structure 
defined as in Lemma \ref{R-cont}.   
Then $A(G,X) \in \mathcal{K}_{iso\mathbb{R}}(\mathsf{gth},\gamma^{m,n}_G)$. \\ 
(b) Any structure of the class 
$\mathcal{K}_{iso\mathbb{R}}(\mathsf{gth},\gamma^{m,n}_G)$ 
is of the form $A(G,X)$ for 
an appropriate isometric action of $G$ 
with hyperbolic continuity moduli $\gamma^{m,n}_G$ 
for the sort $G$. 
\end{theorem} 

{\em Proof.} 
(a) The first part of the theorem is straightforward. 
It uses Lemma \ref{R-cont}. 

(b) Let a continuous structure $M$ belong to 
$\mathcal{K}_{iso\mathbb{R}}(\mathsf{gth}, \gamma^{m,n}_G)$. 
Let $X =\bigcup B_n$, where $B_n$ are sorts of balls of $M$. 
Let $g\in G$ and $x\in B_m$. 
Using axioms 1(a,b) we can find $m< n$ large enough such that 
$\circ_{ms}(g, x, 0) \le m$ for all $s>n$.  

Using axiom 1(c) choose in $B_m$ a sequence 
$u_1 , u_2 ,\ldots, u_i ,\ldots$, 
such that for $n$ and $s$ as above, 
$$
\mathsf{max} (\circ_{ms} (g,x,u_{i+1} ), 
\mathsf{min} (|d(u_i ,u_{i+1}) -\circ_{ms} (g,x,u_i )|< \frac{1}{2^{i+1}} .
$$ 
By completenes of $B_m$ there is $\hat{u} \in B_m$ 
which is the limit of the sequence $u_1 , u_2 ,\ldots, u_i ,\ldots$. 
Using axiom 1(c) and the choice of $u_i$ 
it is easy to see that $\circ_{ms} (g,x,\hat{u} )=0$. 
Define the value $g\circ x$  to be $\hat{u}$. 
By axioms 1(b,c) this value is defined in a unique way. 
This procedure defines a binary function 
$\circ :  G\times X\rightarrow X$. 

By axioms 3(a - c) the function $\circ$ is an action of $G$ on $X$. 
By axioms 2(a) it is easy to see that when $x,y,u \in B_{m}$ 
and $g\circ x = u$ then $\circ_{mn}(g,x,y) = d(u,y)$. 
Using axioms 2(b,c) we have that the action is isometric. 

Since $\circ_{mn}$ has continuity modulus $\gamma^{m,n}_G$ on $G$ 
we easily obtain that the action $\circ$ also has 
the $G$-continuity modulus $\gamma^{m,n}_G$ for the partial function 
$G\times B_m \rightarrow B_n$. 
To see that continuity moduli $\gamma^{m,n}_G$ 
are hyperbolic, note that 
when $x\in B_m$ and $y\in h(B_m )\setminus B_n$, the implication 
$$
\forall g \in G 
(d(1,g)\le \gamma^{m,n}_G (q ) \rightarrow |\mathsf{length} ([x,y]\cap B_n ) 
- \mathsf{length} ([x, g\circ y]\cap B_n )| \le q ).   
$$ 
can be rewritten as 
$$
\forall g \in G 
(d(1,g)\le \gamma^{m,n}_G (q ) \rightarrow |\circ_{mn}(h,h^{-1} (y),x)- \circ_{mn} (gh,h^{-1}(y),x)| \le q ).     
$$ 
We now conclude that $M$ is of the form $A(G,X)$ for 
an appropriate isometric action of $G$ with hyperbolic continuity moduli 
$\gamma^{m,n}_G$. 
$\Box$ 

\bigskip 

We apply this theorem for groups which are non-{\bf F}$\mathbb{R}$. 
Let us start with the following definition which 
modifies non-{\bf F}$\mathbb{R}$ by taking attention to 
the corresponding continuity moduli. 

\begin{definition} 
Let $G$ be a metric group of a bounded diameter 
which is a continuous structure in the language  
$(d,\cdot ,^{-1},1)$.  

We say that $G$ has {\bf property} 
{\bf non-F}$\mathbb{R}$ in the class 
$\mathcal{K}_{iso\mathbb{R}}(\mathsf{gth}, \gamma^{m,n}_G)$
if it has an isometric action on 
an $\mathbb{R}$-tree $X$   
without fixed points and the action 
can be presented as a structure 
$A(G,X) \in \mathcal{K}_{iso\mathbb{R}}(\mathsf{gth}, \gamma^{m,n}_G)$ 
for some point $0\in X$.  
\end{definition} 

If a group is non-{\bf F}$\mathbb{R}$ in 
$\mathcal{K}_{iso\mathbb{R}}(\mathsf{gth}, \gamma^{m,n}_G)$, 
then it is non-{\bf F}$\mathbb{R}$. 
Note that a metric group $(G,d)$ satisfies 
property non-{\bf F}$\mathbb{R}$  in 
$\mathcal{K}_{iso\mathbb{R}}(\mathsf{gth}, \gamma^{m,n}_G)$
if there is a continuous isometric action of 
$G$ on a real tree $X$ and a natural number $s>0$ such that 
the corresponding structure $A(G,X)$ belongs to 
$\mathcal{K}_{iso\mathbb{R}}(\mathsf{gth}, \gamma^{m,n}_G)$ 
and satisfies all statements of the following form:  
$$ 
\mathsf{sup}_{v\in B_{m}} \mathsf{inf}_{g} 
(  {1\over s} \dot{-} \circ_{mn}(g,v,v))= 0  \mbox{ , where } m,n\in \omega \mbox{ and } m<n.    
$$ 
(saying that each element of $B_m$ is moved by some  
element of $G$ by approximately $ {1\over s}$). 
Let us denote it by $\Theta_{m,n,s}$. 
The following proposition shows that this is 
the only reason to be 
non-{\bf F}$\mathbb{R}$. 
It does not use any logic.

\begin{proposition} \label{R-distance} 
If a group $G$ acts on a real tree $X$ by isometries 
without fixed points then there is a natural number $s$ 
such that for any $m$ each element of $B_m$ is moved by some  
element of $G$ by a distance greater than $ {1\over s}$.  
\end{proposition} 

{\em Proof.} 
If $G$ has a hyperbolic element $g$ of hyperbolic length $r$ 
(i.e. $g$ $r$-shifts all points of its axis $L$), 
then it is  easy to see by Lemma \ref{lem_nnest}  
that any element of $X$ is moved by $g$ at the distance $\ge r$. 
Thus $s$ can be chosen so that $ {1\over s}< r$. 

Consider the case when $G$ consists of elliptic elements. 
Since $G$ does not fix any point, by a well-known 
argument $G$ fixes an end (\cite{serre}, Section 6.5, Exercise 2). 
Let $L_0$ be the half-line starting from $0$ which 
represents this end and let $v_1 ,....,v_i ,....$  
be a cofinal $\omega$-sequence in $L_0$ with $d(v_i ,v_{i+1})\ge 1$. 
Then we may assume that $G$ is the union of a strictly 
increasing chain of stabilizers $G_i$ of $v_i$. 

Having $m$ find $j$ with $v_{j-1} \not\in B_m$ (thus $v_j \not\in B_m$). 
Since any arc linking $v_j$ with an element from $B_m$ must 
contain $v_{j-1}$, we see that if $g\in G_j$ fixes a point 
of $B_m$ then it fixes $v_{j-1}$. 
Since $G_j \not= G_{j-1}$ the group $G_j$ contains 
an element $g$ which does not fix any element of $B_m$. 
Since $d(v_j ,v_{j-1}) \ge 1$ any point of $B_m$ can 
be taken by $g$ at a distance greater than 1.   
Thus we define $s=1$. 
$\Box$ 

\bigskip

We now prove logical analyzability of the case of 
non-{\bf F}$\mathbb{R}$ in 
$\mathcal{K}_{iso\mathbb{R}}(\mathsf{gth}, \gamma^{m,n}_G)$.

\begin{theorem} \label{n-FR} 
(a) Every group with property 
non-{\bf F}$\mathbb{R}$ in 
$\mathcal{K}_{iso\mathbb{R}}(\mathsf{gth}, \gamma^{m,n}_G)$   
has an expansion to a structure 
$A(G,X)$ which belongs to the subclass of 
$\mathcal{K}_{iso\mathbb{R}}(\mathsf{gth}, \gamma^{m,n}_G)$ 
axiomatizable by the family 
$\{ \Theta_{s,m,n} : m<n$ , $m,n\in \omega\}$ 
for some fixed $s$. 

(b) 
The class of structures 
$A(G,X)$ witnessing property  
non-{\bf F}$\mathbb{R}$ in 
$\mathcal{K}_{iso\mathbb{R}}(\mathsf{gth}, \gamma^{m,n}_G)$   
is bountiful. 
\end{theorem} 

{\em Proof.} 
Statement (a) follows from Proposition \ref{R-distance}. 

The last statement follows from 
the L\"{o}wenheim-Skolem theorem and 
statement (a).
$\Box$ 

\bigskip 

Note that Theorem \ref{n-FR} states a stronger property 
than just bountifulness of the class of groups 
with property non-{\bf F}$\mathbb{R}$ in 
$\mathcal{K}_{iso\mathbb{R}}(\mathsf{gth}, \gamma^{m,n}_G)$.  
Statement \ref{n-FR}(b) can be considered as a version of 
the corresponding case of Proposition \ref{Bounti}(b) 
in a precise form.

\subsection{Non-FH-actions}

We introduce an axiomatizable class 
$\mathcal{K}_{iso{\bf H}}(\mathsf{gth},\mathsf{ort}^H ,\gamma^{m,n}_G)$ 
of structures of the form $A(G,X)$ 
in the case of Hilbert spaces. 
We will see in Theorem \ref{isoH} that 
the axioms guarantee that all members of 
$\mathcal{K}_{iso{\bf H}}(\mathsf{gth},\mathsf{ort}^H ,\gamma^{m,n}_G)$ 
are obtained in the way described in Lemma \ref{H-cont}. 
After this we show  
logical analyzability of property non-{\bf FH}. 

In some sense the result of this section looks stronger 
than the results of Section 3.2. 
Contrary to the previous case we do not have any additional 
assumptions on continuity moduli of partial functions 
$G\times B_m \rightarrow B_n$. 

\begin{definition} 
Let $\mathcal{K}_{iso{\bf H}}(\mathsf{gth},\mathsf{ort}^H ,\gamma^{m,n}_G)$ 
be a class of continuous metric structures 
which are unions of continuous groups  
$(G, \cdot , ^{-1}, d )$  
together with many-sorted metric structures  
of $n$-balls of real ${\bf H}$   
$$ 
{\bf H} = (\{ B_n\}_{n\in \omega} ,0,\{ I_{mn} \}_{m<n} ,d, \{ \lambda_r \}_{r\in\mathbb{R}}, +,-,\langle \rangle ) 
\mbox{ , where } I_{mn}: B_m\rightarrow B_n ,  
$$ 
and ternary predicates  
$$
\circ_{mn}(g,x,y):G\times B_m \times B_m \rightarrow [0, n+m ]  
\mbox{ , } m, n\in \omega \mbox{ , } \mathsf{ort}^H (m) < n   
$$ 
where $\gamma^{m,n}_G$ is the continuity modulus 
with respect to $G$ and the continuity moduli 
with respect to $B_m$ are equal to $3 \mathsf{id}$. 
The class 
$\mathcal{K}_{iso{\bf H}}(\mathsf{gth},\mathsf{ort}^H ,\gamma^{m,n}_G)$ 
is axiomatizable by the axioms of groups with 
an invariant metric for the sort $G$, 
the axioms of pointed real Hilbert spaces 
and the three groups of axioms of the relation 
$\circ_{mn}$ from Definition \ref{ClassAx} 
where  it is always assumed that $\mathsf{ort}^H (m)<n$ 
and for $x,y\in {\bf H}$ 
the value $\parallel x-y\parallel$ is denoted  
by $d(x,y)$. 
\end{definition}

\bigskip 

The following theorem is a version of 
Theorem \ref{isoR} 
in the case of isometric actions 
on ${\bf H}$. 

\begin{theorem} \label{isoH} 
(a) Assume that a group $G$ acts on ${\bf H} =\bigcup B_n$ 
by affine isometries according to a function $\mathsf{gth}$.  
Let $A(G,{\bf H})$ be the corresponding structure 
defined as in Lemma \ref{H-cont}.   
Assume that the predicates $\circ_{mn}$ for  
$\mathsf{ort}^H (m)<n$, 
have continuity moduli 
$\{ \gamma^{m,n}_G : m<n \in \omega\}$ for $G$-variables 
and $3 \mathsf{id}$ for $B_m$-variables. 
Then 
$A(G,{\bf H}) \in \mathcal{K}_{iso{\bf H}}(\mathsf{gth},\mathsf{ort}^H ,\gamma^{m,n}_G)$. \\ 
(b) Any structure of the class 
$\mathcal{K}_{iso{\bf H}}(\mathsf{gth},\mathsf{ort}^H ,\gamma^{m,n}_G)$ 
is of the form $A(G,{\bf H})$ for 
an appropriate isometric action of $G$ 
on ${\bf H}$. 
\end{theorem} 

{\em Proof.} (a) This part is straightforward. 

(b) The proof of this part repeats the corresponding 
place of the proof of  Theorem \ref{isoR}(b). 
$\Box$ 

\bigskip

We now define an appropriate versions of non-{\bf F}${\bf H}$.  

\begin{definition} \label{def-non-FH}
Let $G$ be a metric group of a bounded diameter 
which is a continuous structure in the language  
$(d,\cdot ,^{-1},1)$.   

(a) We say that $G$ satisfies 
{\bf uniform property non-FH} if there is a function 
$\mathsf{gth}: \mathbb{N} \times [0,1) \rightarrow \mathbb{N}$ 
and there is an affine isometric action $\circ$ of $G$ 
on real ${\bf H}$ without fixed points and  
the action $\circ$ corresponds to $\mathsf{gth}$ and 
has the property that for any $m,n$ with $\mathsf{ort}^H (m)<n$ 
there is a continuity modulus  
for the $G$-variable of the partial function 
$G\times B_m \rightarrow B_n$ induced by $\circ$. 

(b) We say that $G$ has {\bf property} 
{\bf non-F}${\bf H}$ in the class 
$\mathcal{K}_{iso{\bf H}}(\mathsf{gth},\mathsf{ort}^H , \gamma^{m,n}_G)$
if it has an affine isometric action 
on ${\bf H}$ without fixed points and the action 
can be presented as a structure 
$A(G,{\bf H}) \in \mathcal{K}_{iso{\bf H}}(\mathsf{gth},\mathsf{ort}^H , \gamma^{m,n}_G)$.  
\end{definition} 

It is clear that the property of Definition \ref{def-non-FH}(b) 
implies uniform non-{\bf FH} and the latter implies non-{\bf FH}. 
We start logical analysis of these properties with the following remark. 

In the case of 
affine isometric actions on the real Hilbert space 
there is an obvious version of Proposition \ref{R-distance}: 
\begin{quote}
{\em If a group $G$ acts on ${\bf H}$ by isometries 
without fixed points then for any $m$ 
there is a natural number $s$ 
such that each element of $B_m$ is moved by some  
element of $G$ by a distance greater than $ {1\over s}$. } 
\end{quote} 
Indeed since $G$ does not fix any point,  
each orbit of $G$ is unbounded (Proposition 2.2.9 of \cite{BHV}). 
Thus there is $g\in G$ so that $B_m \cap g (B_m) =\emptyset$. 
In particular there is $s\in \mathbb{N}$ such that 
$ {1\over s} \le \parallel g\circ v - v\parallel$ 
for all $v\in B_m$.  

We now see that if a metric group $(G,d)$ satisfies 
property non-{\bf FH} in the class 
$\mathcal{K}_{iso{\bf H}}(\mathsf{gth},\mathsf{ort}^H , \gamma^{m,n}_G)$ 
then there is a $\gamma_G$-continuous affine 
isometric action of $G$ on ${\bf H}$ and 
a function $m\rightarrow s(m)$ such that 
the corresponding structure 
$A(G,{\bf H})$ satisfies all   
statements of the following form:  
$$ 
 \mathsf{inf}_{g} \mathsf{sup}_{v\in B_{m}}
(  {1\over s(m)} \dot{-} \circ_{mn}(g,v,v))= 0  \mbox{ , where } m,n\in \omega \mbox{ and } \mathsf{ort}^H (m)< n.    
$$ 
(saying that there is an element of $G$ which takes 
each element of $B_m$  by approximately $ {1\over s(m)}$). 
Let us denote it by $\Theta^{H,s}_{m,n}$. 

The following theorem gives logical analyzability 
of uniform property non-{\bf FH}.

\begin{theorem} \label{n-FH} 
(a) Every group with property 
non-{\bf FH} in 
$\mathcal{K}_{iso{\bf H}}(\mathsf{gth}, \mathsf{ort}^H , \gamma^{m,n}_G)$   
has an expansion to a structure 
$A(G,{\bf H})$ which belongs to the subclass of 
$\mathcal{K}_{iso{\bf H}}(\mathsf{gth}, \mathsf{ort}^H , \gamma^{m,n}_G)$ 
axiomatizable by the family 
$\{ \Theta^{H,s}_{m,n} : \mathsf{ort}^H (m)<n$ , $m,n\in \omega\}$ 
for some function $s:\mathbb{N} \rightarrow \mathbb{N}$. 

(b) Every group with uniform property non-{\bf FH}
has an expansion to a structure $A(G,{\bf H})$ 
which belongs to the subclass of some 
$\mathcal{K}_{iso{\bf H}}(\mathsf{gth}, \mathsf{ort}^H , \gamma^{m,n}_G)$
axiomatizable by the family 
$\{ \Theta^{H,s}_{m,n} : \mathsf{ort}^H (m)<n$ , $m,n\in \omega\}$ 
for some function $s:\mathbb{N} \rightarrow \mathbb{N}$. 

(c) The class of structures 
$A(G,{\bf H})$ witnessing 
uniform property non-{\bf FH}    
is bountiful. 
\end{theorem}

{\em Proof.} 
Statement (a) follows from the above consequence of 
Proposition 2.2.9 of \cite{BHV}. 
Statement (b) follows from (a), Lemma \ref{H-cont} and Theorem \ref{isoH}.

The last statement follows from the L\"{o}wenheim-Skolem theorem and 
statements (a) and (b).
$\Box$ 


\section{Around non-OB} 

In paragraph A below we uniformize property non-{\bf OB} 
in order to make it logically analyzable. 
Unfortunately when the topology is discrete the class we obtain 
coincides with the class of all infinite groups. 
Since property {\bf OB} is not empty 
in the class of discrete infinite groups 
we consider this class in paragraph B 
using a completely different approach. 
In fact we apply the standard version of $L_{\omega_1 \omega}$. 
These arguments are very close to 
the paper \cite{KM}, 
which was the starting point of 
the subject.

\paragraph{ A. Uniformization.} 
Property {\bf OB} is defined in Definition \ref{def_prop}. 
In this section we study some modifications of {\bf OB} 
in the case of metric groups which are continuous structures.  

It is known (see Section 1.4 of \cite{rosendalN}) 
that for Polish groups property {\bf OB} is equivalent 
to the property that for any open symmetric $V\not=\emptyset$ there is 
a finite set $F\subseteq G$ and a natural number $k$ such that $G=(FV)^k$.  
Thus when $G$ is non-{\bf OB}, there is an non-empty open $V$ such that 
for any finite $F$ and a natural number $k$, $G\not=(FV)^k$. 
Note that for such $F$ and $k$ there is a real number $\varepsilon$ 
such that some $g\in G$ is $\varepsilon$-distant from $(FV)^k$. 
Indeed, otherwise $(FV)^k V$ would cover $G$.  
 
This suggests that the following property can be considered as 
a substitute for the complement to {\bf OB}. 

\begin{definition} 
A metric group $G$ is called {\bf uniformly non-OB} if there is an open 
symmetric $V\not=\emptyset$ such that for any natural numbers $m$ and $k$ 
there is a real number $\varepsilon$ such that for any $m$-element subset 
$F\subset G$ there is $g\in G$ which is $\varepsilon$-distant from $(FV)^k$. 
\end{definition}  

It is clear that 
all infinite discrete groups are uniformly non-{\bf OB}.  
In order to verify that this property is logically analyzable 
we will consider metric groups with two distinguished 
grey subsets of them, i.e. two unary predicates 
denoted by $P$ and $Q$.   
In some sense the definition below states that 
the nullset of $Q$ is the complement of the nullset of $P$. 

\begin{definition} 
{\em Let $\mathcal{K}_{grey}$ be the class of all continuous 
metric structures $\langle G, P,Q \rangle$ such that 
$G$ is a group with a bi-invariant $[0,1]$-metric, 
$P:G\rightarrow [0,1]$ and $Q:G\rightarrow [0,1]$ 
are unary predicates on $G$ with continuity moduli 
$\mathsf{id}$ and the following axioms are 
satisfied:  
$$ 
 \mathsf{sup}_x |P(x)-P(x^{-1})|=0 \mbox{ , } \hspace{2cm} \mathsf{inf}_x |P(x)-1/2| =0,   
$$
$$ 
\mathsf{sup}_x |Q(x)-Q(x^{-1})|=0 \mbox{ , } \hspace{2cm} \mathsf{inf}_x |Q(x)-1/2| =0,    
$$
$$
Q(1)=0 \mbox{ , } \hspace{4cm} \mathsf{sup}_x \mathsf{min} (P(x),Q(x))=0,  
$$
$$
\mathsf{sup}_x  \mathsf{min} ( \varepsilon \dot{-} Q(x), \mathsf{inf}_y (\mathsf{max}(d(x,y)\dot{-} 2\varepsilon , \varepsilon \dot{-}  P(y)))=0,    
$$ 
$$
\hspace{7cm} \mbox{ for all rational }\varepsilon \in [0,\frac{1}{2}].   
$$ 
} 
\end{definition} 

Note that the last axiom implies that any neighbourhood of an element from the nullset of $Q$ 
contains an element with non-zero $P$. 

\begin{definition} 
For any natural $m$ and $k$ and any rational $\varepsilon$ let 
$\theta (m,k,\varepsilon )$ be the following formula: 
$$ 
\mathsf{sup}_{x_1 ...x_m }  \mathsf{inf} _{x}   \mathsf{sup}_{y_1 ...y_k} \mathsf{min}(P(y_1),...,P(y_n), (\varepsilon \dot{-} \mathsf{min}_{w\in W_{m,k}}(d(x,  w))))=0, 
$$ 
$$ 
\mbox{ where } W_{m,k} \mbox{ consists of all words of the form } x_{i_1} y_1 x_{i_2} y_2 ...x_{i_k}y_k . 
$$ 
\end{definition} 

It is worth noting that in the definition the expression  
$\theta (m,k,\varepsilon )$ formally is not a formula. 
On the other hand it is easy to see that 
it can be written by a formula of continuous logic. 
It states that for any finite set $F= \{ x_1 ,...,x_m \}$ 
there is $x$ which is $\varepsilon$-distant 
from $(F \cdot \{ y: P(y) >0\} )^k$. 

The following theorem shows logical analyzability of 
uniform non-{\bf (OB)}. 

\begin{theorem} \label{uniform} 
(a) A metric group $G$ is uniformly non-{\bf OB} if and only if 
there is a finction $s: \mathbb{N} \times \mathbb{N} \rightarrow \mathbb{N}$ 
such that $G$ has an expansion from $\mathcal{K}_{grey}$ which satisfies all 
conditions $\theta (m,k,\frac{1}{s(m,k)} )$. \\ 
(b) The class of uniformly non-{\bf OB} metric groups is bountiful.  
\end{theorem} 
 
{\em Proof.} 
If $G$ is a uniformly non-{\bf OB}-group, then find an open symmetric $V$ 
such that $1\in V$ and for any natural numbers $m$ and $k$ 
there is a real number $\varepsilon$ such that for any $m$-element subset 
$F\subset G$ there is $g\in G$ which is $\varepsilon$-distant from $(FV)^k$. 
We interpret $Q(x)$ by $d(x,V)$ and $P(x)$ by $d(x, G\setminus V)$ 
(possibly normalizing them to satisfy the axioms of $\mathcal{K}_{grey}$). 
Then observe that $\langle G,P,Q\rangle\in \mathcal{K}_{grey}$ and for any 
natural numbers $m$ and $k$ there is a rational number $\varepsilon$ 
so that $\theta(m,k,\varepsilon )$ holds in $(G,P,Q)$.  

Assume that there is a function $s$ such that 
$G$ has an expansion as in statement (a). 
To verify uniform non-{\bf OB} take the complement 
of the nullset of $P(x)$ as an open symmetric subset $V$. 
This proves statement (a).  

Statement (b) follows from (a) and 
the L\"{o}wenheim-Skolem theorem for continuous logic.   
$\Box$ 

\bigskip 

The following definition from \cite{rosendalN} 
gives several versions of property {\bf OB}. 

\begin{definition} 
Let $G$ be a topological group. \\ 
(1) The group $G$  is called bounded if for any open $V$ containing $1$ there is 
a finite set $F\subseteq G$ and a natural number $k>0$  such that $G=FV^k$. \\ 
(2) The group $G$ is Roelcke bounded if for any open $V$ containing $1$ there is 
a finite set $F\subseteq G$ and a natural number $k>0$  such that $G=V^k FV^k$. \\ 
(3) The group $G$ is Roelcke precompact if for any open $V$ containing 
$1$ there is a finite set $F\subseteq G$ such that $G=VFV$. \\ 
(4) The group $G$ has property ${\bf (OB)_k}$ if for any open symmetric 
$V\not=\emptyset$ there is a finite set $F\subseteq G$ such that $G=(FV)^k$. 
\end{definition}

It is worth noting that by Section 1.10 of \cite{rosendalN} 
in the case of $\sigma$-locally compact groups 
(=$\sigma$-compact locally compact)
all these properties coincide with property ${\bf OB}$. 
Applying our method of uniformization 
of the corresponding negations of these properties 
one can obtain bountiful classes of metric groups 
which are:  
\begin{itemize} \label{uniformies} 
\item uniformly non-bounded; 
\item uniformly non-Roelcke bounded;   
\item uniformly non-Roelcke precompact;  
\item uniformly non-${\bf (OB)_k}$.  
\end{itemize}

\paragraph{B. Discrete groups.} 
An abstract group $G$ is {\bf Cayley bounded} if for every generating subset
$U\subset G$ there exists $n\in \omega$ such that every element
of $G$ is a product of $n$ elements of $U\cup U^{-1}\cup\{ 1\}$.
If $G$ is a Polish group then $G$ is {\bf topologically Cayley bounded} 
if for every analytic generating subset $U\subset G$ 
there exists $n\in \omega$ such that every element
of $G$ is a product of $n$ elements of $U\cup U^{-1}\cup\{ 1\}$.
It is proved in \cite{rosendal} that for Polish groups property 
{\bf OB} is equivalent to topological Cayley boundedness together 
with {\bf uncountable topological cofinality}: $G$ is not the union 
of a chain of proper open subgroups. 

Let us consider the abstract (discrete) case. 
Since in this case we do not need continuous logic, 
our considerations become simpler. 

A group is {\bf strongly bounded} if it is Cayley bounded and
cannot be presented as the union of a strictly increasing chain
$\{ H_n :n\in \omega\}$ of proper subgroups
(has {\bf cofinality} $>\omega$).  
It is a discrete version of {\bf OB}, 
which is also called {\bf Bergman's property}, \cite{berg}. 
It is known that strongly bounded groups have property {\bf FA}, 
i.e. any action on a simplicial tree fixes a point. 

As we already know the class of strongly bounded groups is not bountiful. 
The corresponding arguments 
given in Introduction can be also 
applied to property ${\bf FA}$.  

It is shown in \cite{dC}, that strongly bounded groups
have property {\bf FH}.   
It can be also deduced from \cite{dC} that strongly bounded groups 
have property ${\bf F}\mathbb{R}$ that every isometric action of $G$ 
on a real tree has a fixed point (since such a group acting on 
a real tree has a bounded orbit, all the elements are elliptic 
and it remains to apply cofinality $>\omega$). 
It is now clear that the bountiful class of groups having 
free isometric actions on real trees (or on real Hilbert spaces) 
is disjoint from strong boundedness. 
 
In the following Proposition we 
apply the standard version of $L_{\omega_1 \omega}$.  

\begin{proposition} \label{discr} 
The following classes of groups are reducts of axiomatizable 
classes in $L_{\omega_1 \omega}$: \\ 
(1) The complement of the class of strongly bounded groups; \\ 
(2) The class of groups of cofinality $\le \omega$; \\  
(3) The class of groups which are not Cayley bounded; \\ 
(4) The class of groups presented as non-trivial free products with amalgamation 
(or HNN-extensions); \\   
(5) The class of groups having homomorphisms onto $\mathbb{Z}$.  

All these classes are bountiful. 
The class of groups which do not have property {\bf FA} is bountiful too. 
\end{proposition} 

{\em Proof.} 
(1) We use the following characterization of strongly bounded 
groups from \cite{dC}.  
\begin{quote} 
A group is strongly bounded if and only if for every 
presentation of $G$ as $G=\bigcup_{n\in \omega} X_n$ for 
an increasing sequence $X_n$, $n\in \omega$, with 
$\{ 1\} \cup X^{-1}_n \cup X_n \cdot X_n \subset X_{n+1}$ 
there is a number $n$ such that $X_n =G$. 
\end{quote} 
Let us consider the class $\mathcal{K}_{nb}$ of all structures 
$\langle G, X_n \rangle_{n\in\omega}$ with the first-order axioms 
stating that $G$ is a group, $\{ X_n \}$ is a sequence of 
unary predicates on $G$ defining a strictly increasing 
sequence of subsets of $G$ with 
$$
\{ 1\} \cup X^{-1}_n \cup X_n \cdot X_n \subset X_{n+1}
$$ 
and with the following $L_{\omega_1 \omega}$-axiom:  
$$ 
(\forall x) (\bigvee_{n\in\omega} x\in X_n ). 
$$
By the L\"{o}wenheim-Skolem theorem for countable fragments of 
$L_{\omega_1 \omega}$ (\cite{keisler}, p.69)
any subset $C$ of such a structure is contained in an elementary 
submodel of cardinality $|C|$ (the countable fragment which we consider 
is the minimal fragment containing our axioms). 
This proves bountifulness in case (1). \\ 
(2) The case groups of cofinality $\le \omega$ is similar. \\ 
(3) The class of groups which are not Cayley bounded is 
a class of reducts of all groups expanded by an unary predicate  
$\langle G,U\rangle$ with an $L_{\omega_1 \omega}$-axiom stating 
that $U$ generates $G$ and with a system of first-order axioms 
stating that there exists an element of $G$ which is not a 
product of $n$ elements of $U\cup U^{-1} \cup \{ 1 \}$. 
The rest is clear.   \\ 
(4) The class of groups which can be presented as 
non-trivial free products with amalgamation is the class 
of reducts of all groups expanded by two unary predicates  
$\langle G,U_1 ,U_2 \rangle$ with first-order axioms that $U_1$ and $U_2$ 
are subgroups and with $L_{\omega_1 \omega}$-axioms 
stating that $U_1 \cup U_2$ generates $G$ and a word in the 
alphabeth  $U_1 \cup U_2$ is equal to 1 if and only if this 
word follows from the relators of the free product of $U_1$ 
and $U_2$ amalgamated over $U_1 \cap U_2$. 
The rest of (4) is clear.  \\ 
(5) Groups having homomorphisms onto $\mathbb{Z}$ 
can be considered as reducts of structures in the language 
$\langle \cdot ,...U_{-n},...,U_0 ,...,U_{m},...\rangle$, 
where predicates $U_t$ denote preimages of the corresponding integer numbers.

To see that the class of groups without {\bf FA} is 
bountiful, take any infinite $G$ which is not {\bf FA}.  
It is well-known (\cite{serre}, Section 6.1) 
that such a group belongs to the union of 
the classes from statements (2),(4) and (5).  
Thus $G$ has an expansion as in one of the cases (2),(4) or (5). 
Now applying the L\"{o}wenheim-Skolem theorem, for any  
$C\subset G$ we find a subgroup of $G$ of cardinality $|C|$ 
which contains $C$ and does not satisfy {\bf FA}.   
$\Box$

\bigskip

Institute of Mathematics, \parskip0pt

Silesian University of Technology, \parskip0pt

ul. Kaszubska 23, \parskip0pt

44 - 100 Gliwice, \parskip0pt

POLAND. \parskip0pt

{\em iwanowaleksander@gmail.com} 

\end{document}